# COCHAINES QUASI-COMMUTATIVES EN TOPOLOGIE ALGEBRIQUE

**Max Karoubi**

Université Paris 7. UFR de Mathématiques
2, place Jussieu, 75251 Paris Cedex 05
e.mail : karoubi@math.jussieu.fr
http://www.math.jussieu.fr/~karoubi/

Le but de cet article est de construire sur un ensemble simplicial X des “formes différentielles” définies sur un anneau commutatif cohérent[1] k. Elles permettent de définir une nouvelle structure d’algèbre différentielle graduée[2] $\mathcal{D}^*(X)$ dite **quasi-commutative**, qui est quasi-isomorphe à l’algèbre des cochaînes classiques sur X (d’où la terminologie). On montre que cette structure détermine (sous certaines conditions de finitude) le type d’homotopie de X si k = **Z**. En particulier, les opérations de Steenrod sur la cohomologie de X, ainsi que les groupes d’homotopie de X peuvent s’en déduire par des méthodes standard d’algèbre homologique. Notre travail est donc analogue à celui de D. Quillen [24] et D. Sullivan [28] en homotopie rationnelle, où les algèbres différentielles graduées **commutatives** jouent un rôle essentiel. Il est aussi intimement lié à celui de M.A. Mandell sur le type d’homotopie à l’aide des $E_\infty$-algèbres [19], que nous utilisons à la fin de l’article.

Cette structure quasi-commutative sur l’algèbre $\mathcal{D}^*(X)$ enrichit considérablement la théorie classique des cochaînes (notée traditionnellement C*(X)), comme nous comptons le montrer de manière sommaire dans cette introduction. Elle consiste à se donner de manière naturelle, pour tout couple d’espaces X et Y, un sous k-module différentiel gradué $\mathcal{D}^*(X) \overline{\otimes} \mathcal{D}^*(Y)$ de $\mathcal{D}^*(X) \otimes \mathcal{D}^*(Y)$, ainsi qu’un morphisme de modules différentiels gradués (cup-produit)

$$m_{X,Y} : \mathcal{D}^*(X) \overline{\otimes} \mathcal{D}^*(Y) \longrightarrow \mathcal{D}^*(X \times Y),$$

Nous définissons $\mathcal{D}^*(X) \overline{\otimes} \mathcal{D}^*(Y)$ comme le **produit tensoriel réduit** de $\mathcal{D}^*(X)$ et $\mathcal{D}^*(Y)$. Plus généralement, si $X_1, ..., X_n$ sont n espaces, nous définissons leur produit tensoriel réduit

$$\mathcal{D}^*(X_1) \overline{\otimes} \dots \overline{\otimes} \mathcal{D}^*(X_n)$$

comme le sous-module gradué de $\mathcal{D}^*(X_1) \otimes \dots \otimes \mathcal{D}^*(X_n)$ formé de l’intersection des sous-modules $\Gamma_{i,j}$ suivants obtenus à la suite d’une permutation appropriée des facteurs du produit tensoriel

[1] Par exemple un anneau nœthérien.

[2] En fait, dans notre construction $\mathcal{D}^*(X)$ est un k-module plat en raison des hypothèses de cohérence sur k. Pour les applications topologiques, k sera le plus souvent **Z**, $\mathbf{F}_p$ ou **Q**.

$$\Gamma_{i,j} = \ldots \otimes \mathcal{D}^*(X_i)\ \overline{\otimes}\ \mathcal{D}^*(X_j) \otimes \ldots$$

Cette théorie vérifie les propriétés suivantes qui seront précisées dans les deux premiers paragraphes.

1) L'échange des facteurs X et Y induit le diagramme commutatif suivant, où les applications horizontales sont induites par le cup-produit

$$\begin{array}{ccc} \mathcal{D}^*(X)\ \overline{\otimes}\ \mathcal{D}^*(Y) & \rightarrowtail & \mathcal{D}^*(X \times Y) \\ \downarrow & & \downarrow \\ \mathcal{D}^*(Y)\ \overline{\otimes}\ \mathcal{D}^*(X) & \rightarrowtail & \mathcal{D}^*(Y \times X) \end{array}$$

2) L'inclusion naturelle

$$\mathcal{D}^*(X_1)\ \overline{\otimes}\ \ldots\ \overline{\otimes}\ \mathcal{D}^*(X_n) \rightarrowtail \mathcal{D}^*(X_1) \otimes \ldots \otimes \mathcal{D}^*(X_n)$$

est un quasi-isomorphisme.

Il en résulte que toute permutation des facteurs du produit $X_1 \times \ldots \times X_n$ se traduit par la permutation correspondante des facteurs du produit tensoriel réduit

$$\mathcal{D}^*(X_1)\ \overline{\otimes}\ \ldots\ \overline{\otimes}\ \mathcal{D}^*(X_n)$$

de manière compatible avec le cup-produit étendu à n modules

$$\mathcal{D}^*(X_1)\ \overline{\otimes}\ \ldots\ \overline{\otimes}\ \mathcal{D}^*(X_n) \longrightarrow \mathcal{D}^*(X_1 \times \ldots \times X_n)$$

Si les $X_i$ sont réduits à un seul espace X, l'application diagonale de X dans $X^2$ induit une structure d'algèbre différentielle graduée **partielle** sur $\mathcal{D}^*(X)$ dans le sens de Kriz et May : c'est donc une $E_\infty$-algèbre , ce qui permet de lui appliquer la machinerie des opérades [19].

Esquissons maintenant la définition de notre algèbre différentielle graduée $\mathcal{D}^*(X)$ ainsi que celle du produit tensoriel réduit $\overline{\otimes}$ qui joue un rôle essentiel dans cet article.

Comme dans la théorie de Sullivan, il nous faut d'abord trouver l'analogue de l'intervalle [0, 1] en tant qu'algèbre différentielle graduée. En degré 0, celle-ci est l'algèbre $\mathcal{D}^0(x)$ des fonctions

$$f : \mathbf{Z} \longrightarrow k$$

qui sont constantes quand la variable x tend vers $+\infty$ ou $-\infty$ (deux limites indépendantes), le produit étant défini point par point. Le "module des différentielles" $\mathcal{D}^1(x)$ est l'idéal **à gauche** de $\mathcal{D}^0(x)$ formé des combinaisons linéaires finies de fonctions de Dirac (c'est-à-dire partout nulles, sauf en un point où elles sont égales à 1). Si on pose $\overline{f}(x) = f(x + 1)$, la structure de module **à droite** de $\mathcal{D}^1(x)$ s'exprime par la formule suivante

$$\omega.f = \bar{f}.\omega$$

soit $(\omega.f)(x) = f(x + 1)\omega(x)$, avec $f \in \mathcal{D}^0(x)$ et $\omega \in \mathcal{D}^1(x)$.

La différentielle $d : \mathcal{D}^0(x) \longrightarrow \mathcal{D}^1(x)$ est définie par la formule du calcul aux différences

$$(df)(x) = f(x + 1) - f(x)$$

Le k-module $\mathcal{D}^*(x) = \mathcal{D}^0(x) \oplus \mathcal{D}^1(x)$ est alors une algèbre différentielle graduée grâce à la formule de Leibniz qui s'écrit dans sa version non commutative

$$d(f.g) = (df).g + f.(dg)$$

comme le montre un calcul élémentaire classique. Enfin, il est facile de voir que la cohomologie de $\mathcal{D}^*(x)$ est concentrée en degré 0 et est isomorphe à k ("lemme de Poincaré").

L'algèbre commutative $\mathcal{D}^0(x)$ est engendrée librement par la fonction unité et les "fonctions de Heaviside" $Y_x$, $x \in \mathbf{Z}$, définies par $Y_x(n) = 1$ pour $n \leq x$, 0 sinon. Notons que la différentielle de $-Y_x$ est la fonction de Dirac $\omega_x$ égale à 1 si $n = x$ et 0 sinon. Pour la fonction de Heaviside $Y_x$ comme pour la fonction de Dirac $\omega_x$, on dira que $\{x\}$ est leur support singulier.

Le produit tensoriel **réduit** $\mathcal{D}^0(x) \overline{\otimes} \mathcal{D}^0(x)$ est le sous-module libre de $\mathcal{D}^0(x) \otimes \mathcal{D}^0(x)$ de base les produits tensoriels $Y_\alpha \otimes Y_\beta$ avec $\alpha \neq \beta$. On définit de même les produits tensoriels réduits $\mathcal{D}^0(x) \overline{\otimes} \mathcal{D}^1(x)$ (resp. $\mathcal{D}^1(x) \overline{\otimes} \mathcal{D}^0(x)$, resp. $\mathcal{D}^1(x) \overline{\otimes} \mathcal{D}^1(x)$) en considérant les $Y_\alpha \otimes \omega_\beta$ (resp. $\omega_\alpha \otimes Y_\beta$, resp. $\omega_\alpha \otimes \omega_\beta$) avec $\alpha \neq \beta$. Avec cette condition $\alpha \neq \beta$ l'application induite par le cup-produit

$$\mathcal{D}^*(x) \overline{\otimes} \mathcal{D}^*(x) \longrightarrow \mathcal{D}^*(x)$$

est $\mathbf{Z}/2$-équivariante pour la permutation (graduée) des facteurs de $\mathcal{D}^*(x) \overline{\otimes} \mathcal{D}^*(x)$ et pour l'action triviale sur $\mathcal{D}^*(X)$. D'autre part, l'inclusion

$$\mathcal{D}^*(x) \overline{\otimes} \mathcal{D}^*(x) \rightarrowtail \mathcal{D}^*(x) \otimes \mathcal{D}^*(x)$$

est un quasi-isomorphisme. Plus généralement, nous pouvons définir un produit tensoriel réduit de n facteurs $\mathcal{D}^*(x)$ (noté $\mathcal{D}^*(x)^{\overline{\otimes} n}$) avec un quasi-isomorphisme

$$\mathcal{D}^*(x)^{\overline{\otimes} n} \longrightarrow \mathcal{D}^*(x)^{\otimes n}$$

tel que l'application induite par le cup-produit $\mathcal{D}^*(x)^{\overline{\otimes} n} \longrightarrow \mathcal{D}^*(x)$ soit équivariante pour l'action du groupe symétrique sur $\mathcal{D}^*(x)^{\overline{\otimes} n}$ et l'action triviale sur $\mathcal{D}^*(x)$ : on considère des n-produits tensoriels de fonctions de Heaviside et de fonctions de Dirac avec des supports singuliers disjoints deux à deux (cf. 4.3 pour plus de détails).

Dans cette interprétation de l'intervalle [0, 1], $-\infty$ joue le rôle de 0 et $+\infty$ celui de 1, ce qui se traduit par deux augmentations de $\mathcal{D}^*(x)$ vers k obtenues en posant $x = \pm\infty$. De même, il est naturel de définir l'analogue algébrique du (r+1)-cube comme le produit tensoriel gradué

$$\mathcal{D}^*(x_0, \ldots, x_r) = \mathcal{D}^*(x_0) \otimes \mathcal{D}^*(x_1) \otimes \ldots \otimes \mathcal{D}^*(x_r),$$

dont les éléments sont les "formes différentielles" en des variables $x_i$ indépendantes (plus précisément des fonctions adéquates sur $\mathbf{Z}^{r+1}$). De même, le r-simplexe $\Delta_r$ est “défini” par l’équation $x_0 + ... + x_r = +\infty$ (en d’autres termes, une des variables $x_i$ égale à $+\infty$ , en convenant que $t + (+\infty) = +\infty$ pour tout t). De manière précise, $\Omega^*(\Delta_r)$ est l’égalisateur des deux flèches obtenues en rendant égales à $+\infty$ certaines variables

$$\prod_i \Omega^*(x_0, ...., \hat{x}_i, ... , x_r) \rightrightarrows \prod_{i<j} \Omega^*(x_0,..\hat{x}_i,.., \hat{x}_j, ... , x_r)$$

Les opérations face sur $\Omega^*(\Delta_r)$ sont obtenues en posant $x_i = -\infty$ (comparer avec [4]). On notera que pour $i \neq j$, les “variables” $x_i$ et $x_j$, ainsi que leurs différentielles commutent entre elles (au sens gradué). La correspondance $[r] \mapsto \Omega^*(\Delta_r)$ définit une algèbre différérentielle graduée **simpliciale quasi-commutative**, ainsi qu’il est précisé en 2.5 et 4.5. Des versions “bornées” de $\Omega^*(x)$ et de $\Omega^*(\Delta_r)$ sont obtenues en considérant uniquement des fonctions de Heaviside ou de Dirac dont le support singulier est contenu dans un sous-ensemble fixé S de $\mathbf{Z}$ . On notera ces versions bornées ${}_S\Omega^*(x)$ et ${}_S\Omega^*(\Delta_r)$ respectivement.

Soit maintenant $X = X_\bullet$ un ensemble simplicial **fini** (c’est-à-dire ayant un nombre fini de simplexes non dégénérés) L’algèbre des “**formes différentielles”** sur X est définie par l’égalité suivante :

$$\Omega^*(X) = \mathrm{Mor}(X_\bullet, \Omega^*(\Delta_\bullet))$$

où $\mathrm{Mor}(X_\bullet, Y_\bullet)$ désigne en général l’ensemble des morphismes simpliciaux de X vers Y (cf. [3] où de telles “théories cohomologiques” sont décrites, en relation avec des travaux non publiés d’A. Grothendieck). Si X est quelconque (non nécessairement fini), une version légèrement modifiée est présentée dans les § 3 et 4 afin de satisfaire aux axiomes.
Si X est un **complexe** simplicial fini, $\Omega^*(X)$ est aussi l’algèbre obtenue en égalisant les deux morphismes canoniques

$$\prod_\sigma \Omega^*(\sigma) \rightrightarrows \prod_{\tau \subset \sigma} \Omega^*(\tau)$$

$\sigma$ parcourant l’ensemble des simplexes de X (isomorphes à $\Delta_r$ , r = 0, 1...) : la description intuitive de $\Omega^*(X)$ qui en résulte est très proche de celle introduite à l’origine par Sullivan dans le cas rationnel [28].

La cohomologie du complexe $\Omega^*(X)$ est naturellement isomorphe à la cohomologie usuelle de X à coefficients dans k (avec sa structure multiplicative). La structure d’algèbre différentielle graduée quasi-commutative de $\Omega^*(x)$ vue plus haut s’étend en une structure analogue “bisimpliciale” sur les $\Omega^*(\Delta_r)$. Si X et Y sont des ensembles simpliciaux finis, on définit alors le produit tensoriel réduit $\Omega^*(X) \overline{\otimes} \Omega^*(Y)$ comme le k-module des applications **bisimpliciales** de X x Y dans $\Omega^*(\Delta_\bullet) \overline{\otimes} \Omega^*(\Delta_\#)$, ce dernier **produit tensoriel réduit** **bisimplicial** étant défini en 4.5. Là aussi de légères modifications sont nécessaires si X et Y ne sont pas finis.

Les axiomes d'une telle "théorie des cochaînes quasi-commutatives" sont vérifiés dans les paragraphe 3 et 4. Une version légèrement différente (sans doute plus intuitive) dans le contexte des complexes simpliciaux finis est présentée dans le § 5. Des applications sont décrites dans les paragraphes 6 et 7.
Ainsi, dans le § 6, si k = $\mathbf{F}_p$, nous définissons de manière plus conceptuelle les opérations de Steenrod (et plus généralement les cup i-produits) en évitant la machinerie simpliciale, dans le cadre des algèbres différentielles graduées quasi-commutatives.
De manière plus innovante, nous montrons dans le § 7 comment itérer la "bar-construction" en utilisant les produits tensoriels réduits. Ceci permet de décrire un complexe "calculant" les groupes de cohomologie d'espaces de lacets r-itérés d'un espace pointé X (ce qui généralise le travail d'Adams et Hilton pour r = 1). Plus modestement, nous montrons l'existence d'une suite spectrale[3] dont le terme $E_2$ est l'homologie de Hochschild r-itérée et qui converge vers la cohomologie de l'espace des lacets r-itérés de X, noté classiquement $\Omega^r(X)$.
Enfin, comme il a été dit plus haut, nous pouvons associer à cette théorie des cochaînes quasi-commutatives une $E_\infty$-algèbre (dans le sens des opérades), ce qui permet d'utiliser toute la machinerie de [20] pour retrouver le type d'homotopie à partir de cette structure quasi-commutative (si k = **Z**).

Il convient de noter que l'algèbre $\mathcal{D}^*(X)$ a d'autres structures algébriques remarquables qui sont sommairement décrites dans le dernier paragraphe (dans le cadre des **complexes simpliciaux finis**).
Tout d'abord, on montre que l'algèbre $\mathcal{D}^*(X)$ peut être munie d'une structure "tressée", déjà esquissée dans la référence [13]. On obtient ainsi une représentation non triviale du groupe des tresses $\mathbb{B}_n$ sur le module $\mathcal{D}^*(X)^{\otimes n}$. On ignore cependant si cette structure tressée détermine à elle seule le type d'homotopie de X.
Une deuxième structure est un automorphisme T, déduit de la translation usuelle sur **Z**. Celui-ci est homotope à l'identité et vérifie la propriété remarquable suivante : pour tout élément z de $\mathcal{D}^*(X) \otimes \mathcal{D}^*(Y)$, il existe un entier n tel que l'image de l'opérateur $T^n \otimes 1$ appliqué à z appartient au produit tensoriel réduit $\mathcal{D}^*(X) \overline{\otimes} \mathcal{D}^*(Y)$ ; autrement dit, l'opérateur $T \otimes 1$ itéré un certain nombre de fois "fait commuter les cochaînes".
Enfin, on définit de manière naturelle un "produit réduit" $\mathcal{D}^*(X) \,\overline{\mathrm{x}}\, \mathcal{D}^*(Y)$ qui est un **sous-ensemble** du produit ordinaire $\mathcal{D}^*(X) \times \mathcal{D}^*(Y)$ et qui forme un système de générateurs du produit tensoriel réduit $\mathcal{D}^*(X) \overline{\otimes} \mathcal{D}^*(Y)$, considéré comme sous k-module du produit

3 Cette suite spectrale se retrouve aussi (avec d'autres notations) dans un article de Brooke Shipley : convergence of the homology spectral sequence of a cosimplicial space. Amer. J. Math. 118 (1996), N° 1, 179-207.

tensoriel $\Omega^*(X) \otimes \Omega^*(X)$. Ce produit réduit s'étend naturellement à n facteurs.

Enfin, on remarquera que les méthodes de cet article sont suffisamment générales pour s'appliquer à d'autres contextes. Par exemple, si $\mathfrak{F}$ est un faisceau en k-algèbres commutatives, on peut définir de même un "complexe à la de Rham" $\Omega^*(X ; \mathfrak{F})$ qui est une algèbre différentielle graduée quasi-commutative. Il est muni de meilleures propriétés formelles (au niveau du cup-produit) que les différents types de complexes standard associés à un faisceau. D'après la philosophie de Mandell, cette structure supplémentaire permet de donner un sens au "type d'homotopie entier" d'un espace annelé.

On trouvera des variantes et des résumés de cet article dans les références suivantes :

**Remerciements.** Je remercie tout d'abord M. Zisman qui a montré un intérêt constant tout au long des diverses versions de ce texte, tout en apportant des améliorations notables dans la présentation. On peut citer par exemple la remarque 2.8, utilisée souvent dans nos démonstrations. M. Zisman a également écrit un annexe à cet article, sur lequel sont basées les considérations du $ 7 relatives à l'itération de la bar-construction, ainsi que certaines considérations combinatoires du § 4. Je remercie aussi M.A. Mandell pour m'avoir communiqué ses résultats publiés dans [19] et [20] qui sont utilisés de manière essentielle à la fin du § 7 dans la détermination du type d'homotopie à partir de la structure quasi-commutative des algèbres différentielles graduées.

## 1. Généralités sur les algèbres différentielles graduées quasi-commutatives.

L'objet de ce court paragraphe est d'esquisser le cadre algèbrique de cet article. Nous adopterons souvent les abréviations suivantes :

ADG : algèbre différentielle graduée

ADGC : algèbre différentielle graduée **commutative**

ADGQ : algèbre différentielle graduée **quasi-commutative** (cf. 1.1)

Nous fixons un anneau commutatif k dans tout cet article. Pour des raisons de platitude, nous supposerons que k est un anneau cohérent (par exemple un anneau noethérien, comme **Z**, $\mathbf{F}_p$ ou **Q**).

**1.1. DEFINITION.** *Une algèbre différentielle graduée* **quasi-commutative** *est définie par les données suivantes, soumises aux axiomes* α, β *et* γ *explicités plus loin.*

1. *Une algèbre différentielle graduée* A *avec un élément unité* **1** de degré 0.
2. *Un sous* k*-module différentiel gradué* $A^{\overline{\otimes}2}$ *de* $A^{\otimes 2}$, *stable par l'action du groupe* **Z**/2 *opérant naturellement sur* $A^{\otimes 2}$ *et contenant* k.**1**⊗ A *(et donc* A⊗ k.**1**).

Le sous-module $A^{\overline{\otimes}2}$ de $A^{\otimes 2}$ est appelé le **produit tensoriel réduit** ; il est aussi noté $A \overline{\otimes} A$. Par ailleurs, pour i et j appartenant à l'ensemble P = {1, ..., n}, $A_{i,j}$ désigne l'image de $A^{\overline{\otimes}2} \otimes A^{\otimes(n-2)}$ dans $A^{\otimes n}$ grâce une permutation des facteurs du produit tensoriel telle que

$$(1, 2) \mapsto (i, j)$$

Avec cette définition, $A^{\overline{\otimes}n}$ est l'intersection de tous les $A_{i,j}$. Il est aussi noté $A \overline{\otimes} ... \overline{\otimes} A$ (n facteurs)

Voici maintenant les propriétés α, β *et* γ qui caractérisent une ADG quasi-commutative :

**α)** *La restriction de la multiplication* $\mu : A^{\overline{\otimes}2} \longrightarrow A$ *est équivariante, le groupe* **Z**/2 *opérant trivialement sur* A *et par permutation des facteurs sur* $A \overline{\otimes} A$ , *avec les conventions usuelles de signe* (axiome de commutativité)

**β)** *Le* k*-module* $\mu_{12}(A^{\overline{\otimes}3})$ *est inclus*[4] *dans* $A^{\overline{\otimes}2}$, *où* $\mu_{12}$ *désigne le produit sur* $A^{\otimes 3}$ *restreint aux deux premiers facteurs* (axiome d'associativité)[5].

**γ)** *L'inclusion de* $A^{\overline{\otimes}n}$ *dans* $A^{\otimes n}$ *est un quasi-isomorphisme* (axiome de quasi-isomorphisme).

---

[4] Ceci implique par symétrie la même propriété pour $\mu_{23}(A^{\overline{\otimes}3})$

[5] Une ADGQ est donc un cas particulier d'une "ADG partielle", notion définie dans [13] p. 40.

Comme pour les algèbres commutatives, les deux premières propriétés impliquent que toute application f de {1, ..., n} vers {1, ..., p} induit un morphisme fonctoriel

$$f_* : A^{\overline{\otimes} n} \longrightarrow A^{\overline{\otimes} p}.$$

Plus précisément, $f_*$ est induit sur les tenseurs décomposables par la formule suivante

$$f_*(a_1 \otimes ... \otimes a_n) = b_1 \otimes ... \otimes b_p$$

où $b_j$ est le produit de tous les $a_i$ tels que f(i) = j (ce produit est independant de l'ordre des $a_i$ si on se restreint au produit tensoriel réduit).

En particulier, la restriction à $A^{\overline{\otimes} n}$ de la multiplication

$$\mu : A^{\otimes n} \longrightarrow A$$

est équivariante pour l'action naturelle du groupe symétrique $\mathfrak{S}_n$ sur $A^{\overline{\otimes} n}$ et l'action triviale sur A. Cette structure est donc l'analogue algébrique de la notion de Γ-espace introduite par G. Segal [26].

**1.2. Exemples.** Si $A^{\overline{\otimes} 2} = A^{\otimes 2}$, nous retrouvons la notion d'ADG commutative. Comme nous le verrons dans le § 2, l'algèbre A = $\mathcal{D}^*(x)$ considérée dans l'introduction est une ADGQ. Elle servira de base à la construction de $\mathcal{D}^*(X)$ - pour tout ensemble simplicial X - qui est aussi une ADGQ.

**1.3.** Un morphisme entre deux ADG quasi-commutatives A et B est simplement un morphisme d'algèbres différentielles graduées f : A $\longrightarrow$ B tel que $(f \otimes f)(A^{\overline{\otimes} 2}) \subset B^{\overline{\otimes} 2}$ .

**1.4**. Une ADGQ est dite **spéciale** si on se donne en outre un **sous-ensemble** $A_2$ de A x A tel que $A^{\overline{\otimes} 2}$ soit engendré par $A_2$ considéré comme sous-ensemble du produit tensoriel $A^{\otimes 2}$. L'ensemble $A_2$ sera désigné comme le "domaine de commutativité" du produit. On peut aussi définir $A_n$ comme l'intersection de tous les sous-ensembles $A_{i,j}$ qu'on peut considérer dans $A^n$ en privilégiant les facteurs i et j (cf. la définition de $A^{\overline{\otimes} n}$ un peu plus haut). On exige alors que $A^{\overline{\otimes} n}$ est aussi engendré par $A_n$.

## 2. Définition axiomatique d'une théorie des cochaînes quasi-commutatives.

**2.1.** Nous nous plaçons pour simplifier dans la catégorie des ensembles simpliciaux et choisissons par ailleurs un anneau commutatif de base k qui est cohérent (le plus souvent **Z**, $\mathbf{F}_p$ ou **Q** dans les applications). Nous allons tout d'abord reprendre - en les approndissant - les définitions esquissées dans l'introduction et le § 1.

Une **théorie des cochaînes quasi-commutatives** est la donnée d'un foncteur contravariant

$$X \mapsto \mathcal{D}^*(X)$$

de la catégorie des ensembles simpliciaux vers celle des k-modules différentiels gradués **plats**. On se donne aussi un "cup-produit "

$$m_{X,Y} : \mathcal{D}^*(X) \otimes \mathcal{D}^*(Y) \longrightarrow \mathcal{D}^*(X \times Y)$$

ainsi qu'un sous-module différentiel gradué $\mathcal{D}^*(X) \overline{\otimes} \mathcal{D}^*(Y)$ de $\mathcal{D}^*(X) \otimes \mathcal{D}^*(Y)$ (appelé "produit tensoriel réduit", comme dans le § 1), stable par échange des facteurs du produit tensoriel. Ces données doivent satisfaire aux axiomes suivants

**1. Axiome de Mayer-Vietoris**. Un diagramme cocartésien d'ensembles simpliciaux

$$\begin{array}{ccc} X & \longrightarrow & X_1 \\ \downarrow & & \downarrow i_1 \\ X_2 & \xrightarrow{i_2} & X' \end{array}$$

($i_1$ ou $i_2$ étant une cofibration) induit un diagramme cartésien d'algèbres différentielles graduées

$$\begin{array}{ccc} \mathcal{D}^*(X) & \longleftarrow & \mathcal{D}^*(X_1) \\ \uparrow & & \uparrow \\ \mathcal{D}^*(X_2) & \longleftarrow & \mathcal{D}^*(X') \end{array}$$

**2. Axiome d'homotopie et de normalisation**. Si P est un espace réduit à un point, on a $k = \mathcal{D}^*(P)$. Par ailleurs, si I un modèle simplicial de l'intervalle [0, 1], la projection évidente de X x I sur X induit un quasi-isomorphisme entre $\mathcal{D}^*(X)$ et $\mathcal{D}^*(X \times I)$.

D'après Eilenberg et Steenrod, si X est un ensemble simplicial fini, ces deux premiers axiomes impliquent que la cohomologie de $\mathcal{D}^*(X)$ est naturellement isomorphe à la cohomologie ordinaire à coefficients dans k.

**3. Axiome de Künneth**. Le morphisme $m_{X,Y}$ défini plus haut est un quasi-isomorphisme de modules différentiels gradués si X (ou Y) est un ensemble simplicial fini. Par ailleurs, ce morphisme est "associatif" en un sens évident ; enfin il existe un "élément unité" $1 \in \mathcal{D}^0(X)$, définissant un morphisme de k vers $\mathcal{D}^*(X)$, tel qu'on ait un diagramme commutatif :

$$\begin{array}{ccc} \mathcal{D}^*(Y) = k \otimes \mathcal{D}^*(Y) & \longrightarrow & \mathcal{D}^*(X \times Y) \\ \downarrow & & \| \\ \mathcal{D}^*(X) \otimes \mathcal{D}^*(Y) & \longrightarrow & \mathcal{D}^*(X \times Y) \end{array}$$

Dans ce diagramme, la première flèche horizontale est essentiellement induite par la projection de X x Y sur Y. En choisissant X = Y, cet axiome implique que le morphisme composé

$$\mathcal{D}^*(X) \otimes \mathcal{D}^*(X) \xrightarrow{m_{X,X}} \mathcal{D}^*(X \times X) \xrightarrow{\Delta^*} \mathcal{D}^*(X)$$

où Δ est l'application diagonale, définit une structure d'algèbre différentielle graduée (avec élément unité) sur $\mathcal{D}^*(X)$.

**4. Axiome de commutativité**. Le diagramme suivant est commutatif

$$\begin{array}{ccc} \mathcal{D}^*(X) \,\overline{\otimes}\, \mathcal{D}^*(Y) & \longrightarrow & \mathcal{D}^*(X \times Y) \\ R\downarrow & & \sigma\downarrow \\ \mathcal{D}^*(Y) \,\overline{\otimes}\, \mathcal{D}^*(X) & \longrightarrow & \mathcal{D}^*(Y \times X) \end{array}$$

Le morphisme σ est induit par l'échange des facteurs X et Y, tandis que R permute les facteurs du produit tensoriel.

**5. Axiome de quasi-isomorphisme.** Si $X_1, ..., X_n$ sont n espaces, on définit le produit tensoriel réduit

$$\mathcal{D}^*(X_1) \,\overline{\otimes}\, ... \,\overline{\otimes}\, \mathcal{D}^*(X_n)$$

comme le sous-module de $\mathcal{D}^*(X_1) \otimes ...\otimes \mathcal{D}^*(X_n)$ obtenu en considérant l'intersection - pour tous les couples (i, j) - des produits tensoriels réduits $\mathcal{D}^*(X_i) \,\overline{\otimes}\, \mathcal{D}^*(X_j)$ (réordonnés à l'intérieur du produit tensoriel). Ce dernier axiome exprime alors que l'inclusion

$$\mathcal{D}^*(X_1) \,\overline{\otimes}\, ... \,\overline{\otimes}\, \mathcal{D}^*(X_n) \longrightarrow \mathcal{D}^*(X_1) \otimes ...\otimes \mathcal{D}^*(X_n)$$

est un quasi-isomorphisme.

**Remarque 1.** Il résulte de tous ces axiomes que $\mathcal{D}^*(X)$ est une algèbre différentielle graduée quasi-commutative (considérer l'application diagonale de X dans X x X et le morphisme induit sur $\mathcal{D}^*$).

**Remarque 2**. Nous pouvons imposer une condition plus stricte sur la théorie en nous donnant un sous-ensemble $\mathcal{D}^*(X) \,\bar{x}\, \mathcal{D}^*(Y)$ du produit cartésien $\mathcal{D}^*(X) \times \mathcal{D}^*(Y)$ en sorte que $\mathcal{D}^*(X) \,\overline{\otimes}\, \mathcal{D}^*(Y)$ soit engendré par ce "produit réduit" $\mathcal{D}^*(X) \,\bar{x}\, \mathcal{D}^*(Y)$. On imposera aussi que $\mathcal{D}^*(X_1) \,\overline{\otimes}\, ... \,\overline{\otimes}\, \mathcal{D}^*(X_n)$ est engendré par le produit réduit de n facteurs $\mathcal{D}^*(X_1) \,\bar{x}\, ... \,\bar{x}\, \mathcal{D}^*(X_n)$ obtenu en faisant les intersections 2 à 2 des produits réduits à l'intérieur du produit cartésien $\mathcal{D}^*(X_1) \times ... \times \mathcal{D}^*(X_n)$. En suivant la terminologie du § 1, une telle théorie des cochaînes quasi-commutatives est dite spéciale.

**2.2.** Il nous reste maintenant à construire une théorie des cochaînes vérifiant les cinq axiomes (et éventuellement spéciale). Cette construction est détaillée dans ce paragraphe pour le simplexe standard. Le reste de la construction et la vérification des axiomes feront l'objet des deux paragraphes suivants.
Commençons par un examen plus attentif de l'exemple fondamental $\Omega^*(x)$ esquissé dans l'introduction.

**2.3. PROPOSITION** ("lemme de Poincaré" pour la droite affine). *La différentielle*

$$d : \Omega^0(x) \longrightarrow \Omega^1(x)$$

*est surjective et son noyau se réduit aux fonctions constantes.*
*Démonstration.* Si $f \in \mathrm{Ker}(d)$, nous avons $f(x) = f(x+1)$, donc f est une fonction constante. D'autre part, si $\omega(x) \in \Omega^1(x)$, une "primitive" $f(x)$ est la fonction $\sum_{t=-\infty}^{x-1} \omega(t)$ (qui est en fait une somme finie). Il est clair que $f(x+1) - f(x) = \omega(x)$. Notons que la correspondance $\omega \mapsto f$ définit un opérateur d'homotopie I de $\Omega^*(x)$ vers $\Omega^{*-1}(x)$ tel que $dI + Id = 1 - \varepsilon$, où I est défini par 0 sur $\Omega^0(x)$ et où $\varepsilon$ est l'augmentation de $\Omega^0(x) \to k$, obtenue en posant $x = -\infty$ (cf. [14]).

**2.4.** L'algèbre $\Omega^*(x)$ a en fait **deux** augmentations naturelles (poser $x = \pm \infty$). Elles correspondent intuitivement aux deux extrêmités de l'intervalle [0, 1].

**2.5.** Dans l'introduction, nous avons défini l'algèbre des "formes différentielles" sur le (r+1)-cube $K_{r+1}$ comme le produit tensoriel gradué

$$\Omega^*(x_0,\ldots, x_r) = \Omega^*(x_0) \otimes \ldots \otimes \Omega^*(x_r)$$

L'algèbre $\Omega^*(\Delta_r)$ a été aussi définie comme l'égalisateur des deux flèches

$$\prod_i \Omega^*(x_0, \ldots, \hat{x}_i, \ldots, x_r) \rightrightarrows \prod_{i<j} \Omega^*(x_0, \ldots, \hat{x}_i, \ldots, \hat{x}_j, \ldots, x_r)$$

obtenues en posant $x_j$ ou $x_i = +\infty$ .
Pour justifier la notation $\Omega^*(\Delta_r)$, nous devons montrer que la correspondance

$$[r] \mapsto \Omega^*(\Delta_r)$$

définit bien une ADG simpliciale. Pour les opérations face, ceci est évident : il suffit de poser $x_i = -\infty$ (rappelons que $-\infty$ joue le rôle de 0 dans notre interprétation de l'intervalle [0, 1]). Pour les opérations de dégénérescence, nous définissons

$$s_i : \Omega^*(\Delta_r) \longrightarrow \Omega^*(\Delta_{r+1})$$

en ajoutant une variable "consécutive" à la $i^{\text{ème}}$ place. Pour cela, on définit un morphisme $\vartheta$

de $\mathcal{D}^*(t)$ dans $\mathcal{D}^*(t, u)$ qui, composé par les évaluations t ou u = - ∞, redonne l'identité sur $\mathcal{D}^*(t)$ ou $\mathcal{D}^*(u)$. Comme il a été dit dans l'introduction, $\mathcal{D}^0(t)$ est engendré par la fonction unité et les "fonctions de Heaviside" $Y_x$, $x \in \mathbf{Z}$, définies par $Y_x(n) = 1$ pour $n \leq x$ et 0 sinon. Ces fonctions sont liées par les relations $Y_x.Y_z = Y_\alpha$ avec $\alpha = \mathrm{Inf}(x, z)$. Quant aux différentielles $dY_x$, elles sont associées aux fonctions de Heaviside par les relations suivantes :

$$dY_x.Y_z = Y_{z-1}.dY_x$$
$$Y_z.dY_x = dY_x \text{ si } z \geq x$$
$$Y_z.dY_x = 0 \text{ si } z < x$$

[en fait, on a $dY_x = - \omega_x$ où $\omega_x$ est la fonction de Dirac vue comme forme différentielle de degré un]. Puisque les $Y_x$, $dY_x$ (pour $x \in \mathbf{Z}$) et la fonction unité 1 forment une base de $\mathcal{D}^*(t)$, on peut définir

$$\vartheta : \mathcal{D}^*(t) \longrightarrow \mathcal{D}^*(t) \otimes \mathcal{D}^*(u) = \mathcal{D}^*(t, u)$$

par les formules que voici :

$$\vartheta(1) = 1 \otimes 1$$
$$\vartheta(Y_x) = Y_x \otimes Y_x$$
$$\vartheta(dY_x) = dY_x \otimes Y_x + Y_x \otimes dY_x, \text{ soit } \vartheta(\omega_x) = \omega_x \otimes Y_x + Y_x \otimes \omega_x$$

**2.6.** Compte tenu de ces définitions, on vérifie les formules suivantes qui font des $\mathcal{D}^*(\Delta_r)$ une ADG simpliciale (cf. [21]) :

(i) $\partial_i.\partial_j = \partial_{j-1}.\partial_j$ si $i < j$

(ii) $s_i.s_j = s_{j+1}.s_i$ si $i \leq j$

(iii) $\partial_i.s_j = s_{j-1}\partial_i$ si $i < j$

$\partial_j.s_j = \partial_{j+1}.s_j = 1$

$\partial_i.s_j = s_j\partial_{i-1}$ si $i > j + 1$

**2.7. LEMME.** *Le* k*-module* $\mathcal{D}^*(\Delta_r)$ *est projectif.*

*Démonstration.* Nous allons montrer que le morphisme canonique de restriction

$$P : \mathcal{D}^*(x_0, ..., x_r) \longrightarrow \mathcal{D}^*(\Delta_r)$$

admet un inverse à droite (noté i). Pour $f \in \mathcal{D}^*(\Delta_r)$, on pose

$$i(f)(x_0, ..., x_r) = S_1 - S_2 + ... + (-1)^r S_{r+1}$$

où $S_j(x_0, ..., x_r)$ est la somme de tous les $f(x_0, ..., x_r)$ quand j variables $x_\alpha$ sont égales à + ∞ (leur nombre est le coefficient binomial $\binom{r+1}{j}$ ). Par exemple,

$$S_1 = f(+\infty, x_1, ..., x_r) + f(x_0, +\infty, x_2, ..., x_r) + ...$$

.....

$$S_{r+1} = f(+\infty, ..., +\infty)$$

Si nous posons $x_0 = +\infty$, nous retrouvons $f(+\infty, x_1, ..., x_r)$ comme le montre un calcul direct. Le même raisonnement vaut par symétrie pour les autres variables.

**2.8. Remarque** (M. Zisman). Le lemme précédent et sa démonstration fournissent une description plus directe de $\mathcal{D}^*(\Delta_r)$ : c'est le quotient de $\mathcal{D}^*(x_0, ..., x_r)$ par le sous-module $I_r$ des formes différentielles ω telles que $\omega(x_0, ..., x_r) = 0$ si une des variables $x_i$ est égale à$+\infty$ (géométriquement, ceci revient à écrire $\Delta_r$ comme un sous-espace d'un cube de dimension r+1). Ce sous-module a comme base les produits tensoriels

$$f_0 \otimes ... \otimes f_r$$

où les $f_i$ sont des fonctions de Heaviside $Y_x$ (en degré 0) ou des fonctions de Dirac $\omega_x$ (en degré 1). Il en résulte que $\mathcal{D}^*(\Delta_r)$ est en fait un k-module libre.

**2.9. THEOREME.** *La cohomologie de* $\mathcal{D}^*(\Delta_r)$ *est triviale en degré* $> 0$ *et est isomorphe à k en degré 0.*

*Démonstration.* Puisque $\mathcal{D}^*(x_0, ..., x_r)$ est acyclique, il suffit de démontrer que le sous-module $I_r$ précédent est complètement acyclique (toute la cohomologie est triviale). Or $I_r$ n'est autre que le produit tensoriel de r+1 copies du noyau de l'opérateur d'augmentation $\mathcal{D}^*(x) \to k$ , obtenu en posant $x = +\infty$.

**2.10. THEOREME.** *Les groupes d'homotopie du k-module simplicial*

$$r \mapsto \mathcal{D}^*(\Delta_r)$$

*sont triviaux.*

*Démonstration.* Il est facile de voir que ce k-module simplicial est connexe. Par ailleurs, d'après le calcul classique des groupes d'homotopie d'un k-module simplicial, il suffit de vérifier la propriété d'extension suivante :

Soit $\omega \in \mathcal{D}^*(\Delta_{r-1})$ avec $\partial_i(\omega) = 0$ pour i = 0, ..., r-1 (les $\partial_i$ étant les opérateurs face usuels). Il existe alors $\tilde{\omega} \in \mathcal{D}^*(\Delta_r)$ tel que $\partial_i(\tilde{\omega}) = 0$ pour i = 0, ..., r-1 et $\partial_r(\tilde{\omega}) = \omega$. Pour cela, on pose $\tilde{\omega} = Y(x_r)\,\omega$, où $Y = Y_0$ est la fonction de Heaviside ayant 0 comme support singulier ; il est clair que $\tilde{\omega}$ vérifie les conditions requises.

**2.11.** Un examen plus attentif des définitions précédentes montrent qu'on peut les généraliser pour un sous-ensemble quelconque S de **Z**, en lui associant la sous-algèbre différentielle graduée ${}_S\mathcal{D}^*(x)$ de $\mathcal{D}^*(x)$ engendrée, pour $s \in S$, par les fonctions de Heaviside $Y_s$ , les fonctions de Dirac $\omega_s$ et la fonction unité 1. Les relations entre les Y et les ω sont alors les suivantes

$Y_u Y_t = Y_{\inf(u, t)}$
$Y_u \omega_t = 0$ si $u < t$
$Y_u \omega_t = \omega_t$ si $u \geq t$
$\omega_t Y_u = 0$ si $u < t + 1$

$\omega_t Y_u = \omega_t$ si $u \geq t + 1$

$\omega_u \omega_t = 0$

$dY_u = - \omega_u$

(nous gardons le signe - en raison de l'origine "historique" de la définition).

On notera le cas particulier où S est réduit à un élément. L'algèbre ${}_S\mathcal{D}^*(x)$ est alors engendrée par $Y = Y_0$ et $\omega = \omega_0$ avec les relations $Y^2 = Y$, $Y\omega = \omega$ et et $\omega Y = 0$.

Dans ce cas particulier, on a la relation $dY.Y = (1 - Y).dY$, conséquence de l'identité $Y^2 = Y$, avec en plus $Y.dY = dY$ (d'où $dY.Y = 0$ avec $dY = -\omega$).

La définition que nous avons donnée de l'ADG simpliciale $\mathcal{D}^*(\Delta_r)$ et ses propriétés s'étendent sans problème au cas des formes différentielles à support S quelconque. Nous noterons ${}_S\mathcal{D}^*(\Delta_r)$ cette extension de la définition.

## 3. Construction des algèbres $\mathcal{D}^*(X)$ et ${}_S\mathcal{D}^*(X)$ .

**3.1.** En suivant Henri Cartan (cf. [3] ou [9]), nous pouvons utiliser les théorèmes 2.9 et 2.10 (et plus généralement leurs analogues pour la théorie ${}_S\mathcal{D}^*$) pour définir la cohomologie à coefficients dans k. Si $X = (X_\bullet)$ est un ensemble simplicial quelconque, la méthode de Cartan nous permet de construire[6] une ADG de "formes différentielles simpliciales" ${}_S\underline{\mathcal{D}}^*(X)$ : ses éléments sont les applications simpliciales de $X_\bullet$ dans l'ADG simpliciale ${}_S\mathcal{D}^* = {}_S\mathcal{D}^*(\Delta_\bullet)$. Le fait qu'on obtient bien une théorie cohomologique résulte des axiomes de [3], établis précisément en 2.9 et 2.10 (vrais aussi pour les formes à support quelconque).

**3.2.** On peut interpréter le théorème de Cartan dans le cadre de complexes de chaînes, grâce à la correspondance de Dold-Kan [29]. De manière précise, si $C_*$ et $S_*$ sont deux k-modules simpliciaux, $\mathrm{Hom}(C_*, S_*)$ s'identifie au k-module $\mathrm{Hom}(\overline{C}_*, \overline{S}_*)$ des morphismes entre les complexes normalisés $\overline{C}_*$ et $\overline{S}_*$ (qui sont des complexes de chaînes). On obtient alors le théorème suivant (à comparer avec 3.8 un peu plus loin).

**3.3. THEOREME.** *Supposons que le complexe* $\overline{S}_*$ *ait une homologie triviale et que* $S_n$ *soit un* k*-module projectif pour tout* n. *Une suite exacte de k-modules simpliciaux* :

$$0 \longrightarrow C'_* \longrightarrow C_* \longrightarrow C''_* \longrightarrow 0$$

*induit alors une suite exacte*

$$0 \longrightarrow \mathrm{Hom}(\overline{C}''_*, \overline{S}_*) \longrightarrow \mathrm{Hom}(\overline{C}_*, \overline{S}_*) \longrightarrow \mathrm{Hom}(\overline{C}'_*, \overline{S}_*) \longrightarrow 0$$

---

[6] Nous utilisons la notation $\underline{\mathcal{D}^*}$ au lieu de $\mathcal{D}^*$, car cette définition n'est pas exactement celle que nous choisirons par la suite. Cependant, si X est fini, les deux définitions seront les mêmes.

*Démonstration*. D'une manière générale, considérons le complexe en homomorphismes $\mathcal{HOM}(\overline{C}_*, \overline{S}_*)$, où un terme de degré (homologique) r est donné par une famille f d'applications de degré -r. La différentielle de f est alors donnée par la formule usuelle

$$(df)(x) = d(f(x)) + (-1)^r f(dx)$$

On va montrer que l'homologie du complexe $\mathcal{HOM}$ est réduite à 0 en tout degré. En effet, puisque $S_n$ est un module projectif, il en est de même de $\overline{S}_n$. En outre, le complexe $\overline{S}_*$ admet un opérateur d'homotopie évident obtenu en écrivant le k-module des cycles en facteur direct dans le complexe total. Il en résulte par composition un opérateur d'homotopie pour le complexe $\mathcal{HOM}(\overline{C}_*, \overline{S}_*)$ qui est donc totalement acyclique.

Le théorème se ramène à démontrer l'exactitude de la suite suivante au niveau des cycles de degré 0 du complexe $\mathcal{HOM}$ :

$$0 \longrightarrow Z_0(\mathcal{HOM}(\overline{C}''_*, \overline{S}_*)) \longrightarrow Z_0(\mathcal{HOM}(\overline{C}_*, \overline{S}_*)) \longrightarrow Z_0(\mathcal{HOM}(\overline{C}'_*, \overline{S}_*)) \longrightarrow 0$$

Seule la sujectivité n'est pas évidente. Elle résulte de la nullité du groupe d'homologie $H_0(\mathcal{HOM}(\overline{C}'_*, \overline{S}_*))$ démontrée précédemment.

**3.4.** Appliquons ces généralités à un diagramme cocartésien (C) d'ensembles simpliciaux

$$\begin{array}{ccc} X & \longrightarrow & X_1 \\ \downarrow & & \downarrow i_1 \\ X_2 & \xrightarrow{i_2} & X' \end{array}$$

($i_1$ ou $i_2$ étant une cofibration). On en déduit une suite exacte des k-modules simpliciaux de chaînes

$$0 \longrightarrow C_*(X) \longrightarrow C_*(X_1) \oplus C_*(X_2) \longrightarrow C_*(X') \longrightarrow 0$$

ainsi qu'un diagramme cartésien d'algèbres différentielles graduées d'après le théorème 3.3

$$\begin{array}{ccc} s\underline{\mathcal{D}}^*(X) & \longleftarrow & s\underline{\mathcal{D}}^*(X_1) \\ \uparrow & & \uparrow \\ s\underline{\mathcal{D}}^*(X_2) & \longleftarrow & s\underline{\mathcal{D}}^*(X') \end{array}$$

L'axiome de Mayer-Vietoris s'en déduit.

**3.5.** Les axiomes d'homotopie et de normalisation résultent immédiatement des considérations précédentes. Cependant, la vérification des axiomes 4 et 5 décrits en 2.1 nécessite une modification de la définition de Cartan qu'il nous faut maintenant préciser.

**3.6.** Nous allons d'abord introduire une notion générale de "coproduit tensoriel" entre k-modules simpliciaux et cosimpliciaux. De manière précise, soient $C^*$ un k-module cosimplicial et $S_*$ un module simplicial. Soient $\overline{C}^*$ et $\overline{S}_*$ les modules normalisés de $C^*$ et $S_*$ respectivement. Ce sont des complexes de cochaînes ou de chaînes (d'après la correspondance de Dold-Kan de nouveau) qui sont facteurs directs canoniques de $C^*$ et $S_*$.

Nous définirons le **produit restreint** $\prod_n{}' C^n \otimes S_n$ comme l'ensemble des suites $\omega = (\omega_n)$ dans le produit ordinaire $\prod_n C^n \otimes S_n$ telles que leurs images dans le produit des normalisés $\prod_n \overline{C}^n \otimes \overline{S}_n$ appartiennent en fait à la **somme** directe $\bigoplus_n \overline{C}^n \otimes \overline{S}_n$ .

Le **coproduit tensoriel** des modules $C^*$ et $S_*$, - que nous noterons $C^* \Delta S_*$ - est alors l'ensemble des suites $\theta$ dans le produit restreint $\prod_n{}' C^n \otimes S_n$ qui sont dans le noyau de $(u_* \otimes 1 - 1 \otimes u^*)$. Ici $u_*$ et $u^*$ sont les morphismes standard $u_* : S_n \to S_p$ et $u^* : C^p \to C^n$ induits par une application non décroissante $u : [p] \longrightarrow [n]$.

Si $S_n$ est un k-module projectif de type fini pour tout n et si $S^*$ désigne le k-module cosimplicial dual, on a un isomorphisme canonique

$$C^* \Delta S_* \cong \mathrm{Hom}_f(S^*, C^*)$$

Ici $\mathrm{Hom}_f$ désigne le k-module des morphismes cosimpliciaux $f_n : S^n \to C^n$ tels que $f_n(s)$ soit une **combinaison linéaire** d'éléments dégénérés si n est assez grand. Avec les mêmes hypothèses de finitude sur $S_*$ et si $C^*$ est le dual d'un module simplicial $C_*$, on a aussi des isomorphismes canoniques

$$C^* \Delta S_* \cong \mathrm{Hom}_f(C_*, S_*) \cong \mathrm{Hom}_f(S^*, C^*) \cong S^* \Delta C_*$$

**3.7.** Nous définissons de même le coproduit tensoriel $\overline{C}^* \Delta \overline{S}_*$ d'un complexe de cochaînes et de chaînes en considérant dans la **somme directe** $\bigoplus_n \overline{C}^n \otimes \overline{S}_n$, le noyau de $d \otimes 1 - 1 \otimes d'$ en utilisant des notations analogues aux précédentes. Avec les hypothèses de finitude sur $S_*$ mentionnées précédemment et grâce à la correspondance de Dold Kan de nouveau, on a les isomorphismes suivants ( $C^*$ étant le le dual d'un k-module simplicial $C_*$)

$$C^* \Delta S_* \cong \mathrm{Hom}_f(C_*, S_*) \cong \mathrm{Hom}_f(\overline{C}_*, \overline{S}_*) \cong \overline{C}^* \Delta \overline{S}_*$$

Les symboles $\mathrm{Hom}_f$ représentent ici indifféremment des ensembles de morphismes "finis" entre k-modules simpliciaux ou entre complexes de k-modules. L'énoncé suivant est alors une reformulation "finie" du théorème 3.3 et se démontre de manière analogue.

**3.8. THEOREME.** *Supposons que le complexe* $\overline{S}_*$ *ait une homologie triviale et que* $S_n$ *soit un* k*-module projectif de type fini . Une suite exacte de k-modules cosimpliciaux*

$$0 \longrightarrow C'^* \longrightarrow C^* \longrightarrow C''^* \longrightarrow 0$$

*induit alors une suite exacte des coproduits tensoriels associés*

$$0 \longrightarrow C'^* \,\Delta\, S_* \longrightarrow C^* \,\Delta\, S_* \longrightarrow C''^* \,\Delta\, S_* \longrightarrow 0$$

*De même, si* $S^*$ *est un module cosimplicial qui est projectif de type fini en chaque degré et a une cohomologie triviale, une suite exacte de* k*-modules simpliciaux*

$$0 \longrightarrow C'_* \longrightarrow C_* \longrightarrow C''_* \longrightarrow 0$$

*induit une suite exacte des coproduits tensoriels associés*

$$0 \longrightarrow S^* \,\Delta\, C'_* \longrightarrow S^* \,\Delta\, C_* \longrightarrow S^* \,\Delta\, C''_* \longrightarrow 0$$

**3.9.** Si X est fini, les généralités précédentes permettent d'interpréter ${}_S\underline{\Omega}^*(X)$ $= \mathrm{Hom}(C(X_\#), {}_S\Omega^*(\Delta_\#))$ comme le coproduit tensoriel suivant

$$C^\#(X) \,\Delta\, {}_S\Omega^*(\Delta_\#) \cong \mathrm{Hom}_f(C(X_\#), {}_S\Omega^*(\Delta_\#))$$

Si X est infini et S fini, nous choisirons comme "bonne" définition de ${}_S\Omega^*(X)$ le coproduit tensoriel

$$C^\#(X) \,\Delta\, {}_S\Omega^*(\Delta_\#)$$

Par définition, le k-module $\Omega^*(X)$ est alors la limite inductive des ${}_S\Omega^*(X)$ lorsque S parcourt les sous-ensembles finis de **Z**.

Plus généralement, si S est fini et si $M^\#$ est un k-module cosimplicial plat, on définira[7] ${}_S\Omega^*(M)$ comme le coproduit tensoriel

$${}_S\Omega^*(M) = M^\# \,\Delta\, {}_S\Omega^*(\Delta_\#)$$

On a aussi ${}_S\Omega^*(M) \cong \mathrm{Hom}_f({}_S\Omega_*(\Delta^\#), M^\#)$, où ${}_S\Omega_*(\Delta^\#)$ est le k-module dual de ${}_S\Omega^*(\Delta_\#)$.

Comme ci-dessus, on définit enfin $\Omega^*(M)$ comme la limite inductive des ${}_S\Omega^*(M)$, S parcourant les sous-ensembles finis de **Z**.

On va maintenant montrer que la cohomologie de ${}_S\Omega^*(M)$ et de $\Omega^*(M)$ définit une "théorie de la cohomologie" par rapport à M en un sens que nous allons préciser.

**3.10.** Puisque $\pi_p(Z^p(\Delta_*)) \cong k$ ($Z^p$ désigne ici comme dans [9] le k-module simplicial des "formes différentielles" fermées dans ${}_S\Omega^*(\Delta_*)$), nous pouvons choisir un représentant $\chi_p \in Z^p(\Delta_p)$ qui est 0 sur toutes les faces. Comme il est bien connu (cf. [10] § 3 par exemple), ces formes $\chi_p$ peuvent être déterminées par récurrence sur p, en commençant par le choix évident

[7] avec un léger abus de notation.

de $\chi_0$. De manière générale, nous écrivons $\chi_p$ comme la restriction à la 0-face $\Delta_p$ d'une forme $\omega_{p+1}$ appartenant à ${}_S\mathcal{A}^p(\Delta_{p+1})$, nulle sur toutes les autres faces (c'est en fait la définition de la normalisation ${}_S\overline{\mathcal{A}}(\Delta_{p+1})$). Au cran suivant, nous choisissons $\chi_{p+1}$ comme $d\omega_{p+1}$.
Nous définissons finalement un morphisme

$$\theta_p : \quad \overline{M}^p \longrightarrow \overline{M}^* \,\Delta\, {}_S\overline{\mathcal{A}}^p(\Delta_*)$$

par la formule que voici

$$\theta_p(m) = m \otimes \chi_p + (-1)^{p+1}\, dm \otimes \omega_{p+1}$$

(on rappelle que les éléments de $\overline{M}^* \,\Delta\, {}_S\overline{\mathcal{A}}^p(\Delta_*)$ sont décrits comme des cycles du produit tensoriel des complexes $\overline{M}^* \otimes {}_S\overline{\mathcal{A}}^p(\Delta_*)$).
Puisque le diagramme

$$\begin{array}{ccc} \overline{M}^p & \longrightarrow & \overline{M}^* \,\Delta\, {}_S\overline{\mathcal{A}}^p(\Delta_*) \\ \downarrow & & \downarrow \\ \overline{M}^{p+1} & \longrightarrow & \overline{M}^* \,\Delta\, {}_S\overline{\mathcal{A}}^{p+1}(\Delta_*) \end{array}$$

est commutatif, les $\theta_p$ définissent bien un morphisme de complexes de cochaînes.

**3.11. THEOREME**. *Supposons que* S *soit fini. Le morphisme* $\theta$ *défini ci-dessus induit un quasi-isomorphisme entre les complexes de cochaînes* $\overline{M}^\natural$ *et* $\overline{M}^* \,\Delta\, {}_S\overline{\mathcal{A}}^\natural(\Delta_*)$.

*Démonstration.* Le complexe $\overline{M}^*$ est de manière évidente une limite inductive de complexes de longueur finie et la limite inductive commute avec le coproduit tensoriel des complexes (car nous considérons des **sommes directes**). D'après 3.8, il suffit donc de démontrer le théorème lorsque $\overline{M}^*$ est concentré en un seul degré n. Le complexe $\overline{M}^* \,\Delta\, {}_S\overline{\mathcal{A}}^\natural(\Delta_*)$ s'identifie alors au produit tensoriel $\overline{M}^n \otimes {}_S\overline{\overline{\mathcal{A}}}^\natural(\Delta_n)$, où ${}_S\overline{\overline{\mathcal{A}}}^\natural(\Delta_n)$ est l'espace des "formes différentielles" sur $\Delta_n$ qui s'annulent sur **toutes** les faces. Puisque $\overline{M}^n$ est plat, ce complexe a une cohomologie égale à $\overline{M}^n \otimes \tilde{H}^n(\Sigma^n)$, où $\Sigma^n$ est la sphère de dimension n (vue comme quotient de $\Delta_n$ par sa frontière) et $\tilde{H}^n(\Sigma^n)$ sa cohomologie réduite. Le théorème 3.11 en résulte aussitôt.
Le théorème suivant est une conséquence immédiate des considérations précédentes.

**3.12. THEOREME.** *Soit* $C^*(X)$ *le complexe de cochaînes associé à l'ensemble simplicial* X *et soit* ${}_S\mathcal{A}^\natural(X)$ *le complexe de* k*-modules* $C^*(X) \,\Delta\, {}_S\mathcal{A}^\natural(\Delta_*)$. *Nous avons alors un diagramme commutatif naturel de quasi-isomorphismes*

$$\begin{array}{ccc} {}_S\mathcal{A}^\natural(X) = \; C^*(X) \,\Delta\, {}_S\mathcal{A}^\natural(\Delta_*) & \longrightarrow & \mathrm{Hom}(X_*, {}_S\mathcal{A}^\natural(\Delta_*)) = {}_S\underline{\mathcal{A}}^\natural(X) \\ \nwarrow & & \nearrow \\ & C^\natural(X) & \end{array}$$

*En particulier, pour* $S \subset T$ *l'inclusion de* ${}_S\mathcal{A}^\natural(X)$ *dans* ${}_T\mathcal{A}^\natural(X)$ *est un quasi-isomorphisme.*

**3.13. Remarque.** Le morphisme canonique $C^{\natural}(X) \longrightarrow {}_S\Omega^{\natural}(X)$ défini ci-dessus est un morphisme de modules différentiels gradués et <u>non</u> un morphisme d'ADG. En suivant la méthode décrite dans [9], on peut cependant trouver une suite en zigzag de quasi-isomorphismes d'ADG, mais nous n'en aurons pas besoin ici.

## 4. Cup-produit et produit tensoriel réduit.

**4.1.** Si S est fini, le cup-produit

$$m_{X,Y} : {}_S\Omega^n(X) \otimes {}_S\Omega^m(Y) \longrightarrow {}_S\Omega^{n+m}(X \times Y)$$

est obtenu en composant le morphisme évident

$$\phi : \mathrm{Mor}_f(X_{\natural}, {}_S\Omega^n(\Delta_{\natural})) \otimes \mathrm{Mor}_f(Y_{\natural}, {}_S\Omega^m\Delta_{\natural})) \to \mathrm{Mor}_f(X_{\natural} \times Y_{\natural}, {}_S\Omega^n(\Delta_{\natural}) \otimes {}_S\Omega^m(\Delta_{\natural}))$$

par l'application produit (non commutative)

$${}_S\Omega^n(\Delta_{\natural}) \otimes {}_S\Omega^m(\Delta_{\natural}) \longrightarrow {}_S\Omega^{n+m}(\Delta_{\natural})$$

(cf. [3]). On peut remarquer que la première flèche se factorise par le k-module des morphismes **bisimpliciaux**

$$\mathrm{Mor}_f(X_{\natural} \times Y_{\#}, {}_S\Omega^n(\Delta_{\natural}) \otimes {}_S\Omega^m(\Delta_{\#}))$$

qu'on restreint ensuite à la diagonale. Pour vérifier que la flèche $m_{X,Y}$ est un quasi-isomorphisme si Y est fini (par exemple), on peut raisonner par récurrence sur le nombre de cellules de Y (en utilisant l'axiome de Mayer-Vietoris). Il suffit donc de vérifier l'assertion pour $Y = \Delta_r$, ce qui est clair puisque ${}_S\Omega^*(X \times \Delta_r)$ est quasi-isomorphe à ${}_S\Omega^*(X)$.
Si S est infini, on convient de considérer la limite inductive des ${}_T\Omega^*$ lorsque T parcourt les sous-ensembles finis de S.

**4.2.** Nous allons maintenant vérifier les axiomes de quasi-isomorphisme et de commutativité qui sont les plus délicats de la théorie, car ils vont nécessiter des hypothèses sur S. Nous commençons par les vérifications les plus élémentaires

**4.3. PROPOSITION.** *L'algèbre différentielle graduée* $\Omega^*(x)$ *considérée dans l'introduction est quasi-commutative* (avec les définitions du § 1). *Plus généralement le morphisme structurel* ${}_S\Omega^*(x)^{\overline{\otimes} n} \longrightarrow {}_S\Omega^*(x)^{\otimes n}$ *est un quasi-isomorphisme pour tout sous-ensemble* S *de* **Z**.
*Démonstration.* La seule vérification non triviale est le quasi-isomorphisme

$${}_S\Omega^*(x)^{\overline{\otimes} n} \longrightarrow {}_S\Omega^*(x)^{\otimes n}$$

ce qui revient à démontrer le lemme de Poincaré pour le module ${}_S\Omega^*(x)^{\overline{\otimes} n}$. Nous rappelons que $M = {}_S\mathcal{D}^*(x)^{\overline{\otimes} n}$ est engendré par les produits tensoriels de 1, de fonctions de Heaviside $Y_\alpha$ et de Dirac $\omega_\alpha$ telles que leurs supports singuliers[8] soient disjoints deux à deux et contenus dans S. Notons $M_r$ , $r < n$, le sous-module de M engendré par les produits tensoriels $\omega = 1 \otimes ... \otimes 1 \otimes f \otimes g ...$ , où f est une fonction de Heaviside ou une fonction de Dirac et où le nombre de 1 est égal à r. Il est clair que M est la somme directe de k et des $M_r$ et que cette décomposition en somme directe est compatible avec la différentielle. Définissons maintenant un opérateur d'homotopie $K : M_r \longrightarrow M_r$ en associant à

$$1 \otimes ... \otimes 1 \otimes \omega_x \otimes g \otimes ... \text{ le produit tensoriel } 1 \otimes ... \otimes 1 \otimes (-Y_x) \otimes g ....$$

On vérifie facilement que $dK + Kd = 1$, ce qui implique la trivialité de la cohomologie de tous les $M_r$ ; celle de M est donc isomorphe à k. Par conséquent, la cohomologie de $\mathcal{D}^*(x)^{\overline{\otimes} n}$ est isomorphe à celle de $\mathcal{D}^*(x)^{\otimes n}$ , soit à k en degré 0 et 0 sinon, par la formule de Künneth.

Une autre démonstration du quasi-isomorphisme précédent consiste à introduire un opérateur de projection compatible avec la différentielle

$$P : \mathcal{D}^*(x)^{\otimes n} \longrightarrow \mathcal{D}^*(x)^{\overline{\otimes} n}$$

Il est défini de la manière suivante (où les $f_i$ sont 1, des fonctions de Heaviside ou des fonctions de Dirac) :

$P(f_1 \otimes .. \otimes f_n) = f_1 \otimes .. \otimes f_n$ si les supports singuliers des $f_i$ et $f_j$ sont disjoints pour $i \neq j$ et $P(f_1 \otimes .. \otimes f_n) = 0$ sinon.

**4.3.1. Remarque.** En fait, l'algèbre $\mathcal{D}^*(x)$ est spéciale : le **produit réduit** $\mathcal{D}^*(x) \,\overline{\mathrm{x}}\, \mathcal{D}^*(x)$ peut être défini comme l'ensemble formé des couples (f, g) où f et g sont des fonctions de Dirac ou de Heaviside de supports singuliers disjoints, ou bien des couples (f, g) où l'une des fonctions est égale à 1.

**4.4.** La proposition précédente se généralise au produit tensoriel de plusieurs facteurs. A cet effet, posons ${}_S\mathcal{D}^*(x_0,..., x_r) = {}_S\mathcal{D}^*(x_0) \otimes ... \otimes {}_S\mathcal{D}^*(x_r)$. On peut munir alors ${}_S\mathcal{D}^*(x_0,..., x_r)$ d'une structure d'ADG quasi-commutative en définissant ${}_S\mathcal{D}^*(x_0,..., x_r)^{\overline{\otimes} n}$ comme le sous-module de ${}_S\mathcal{D}^*(x_0,..., x_r)^{\otimes n}$ de base les n-produits tensoriels de tenseurs décomposables $\omega_\alpha$ à supports singuliers[9] disjoints 2 à 2. L'opérateur de projection défini plus haut se généralise aisément dans ce contexte pour montrer que ce sous-module a bien la bonne cohomologie.

Plus généralement, on définit un produit tensoriel réduit

---

[8] Il convient de noter que ce module se réduit à k si le cardinal de S est plus petit que n.

[9] Le support singulier d'un tenseur décomposable $f_0 \otimes f_1 \otimes ... \otimes f_n$ de $\mathcal{D}^*(x_0) \otimes \mathcal{D}^*(x_1) \otimes ... \otimes \mathcal{D}^*(x_n)$ est la réunion des supports singuliers des $f_i$.

$$ {}_{S_1}\mathcal{D}^*(X_1)\ \overline{\otimes}\ {}_{S_2}\mathcal{D}^*(X_2)\ \overline{\otimes}\ ...\ \overline{\otimes}\ {}_{S_m}\mathcal{D}^*(X_m) $$

où les $X_i$ sont des groupes de variables choisies parmi $x_0, ..., x_r$ : il suffit de considérer des combinaisons linéaires de produits tensoriels décomposables en fonctions de Heaviside et de Dirac, tels que les m supports singuliers soient disjoints 2 à 2 et contenus dans $S_1, .. ,$ $S_m$respectivement. Enfin, on définit un produit réduit

$$ {}_{S_1}\mathcal{D}^*(X_1)\ \overline{x}\ {}_{S_2}\mathcal{D}^*(X_2)\ \overline{x}\ ...\ \overline{x}\ {}_{S_m}\mathcal{D}^*(X_m) $$

de manière tout à fait analogue : ce produit réduit engendre le produit tensoriel réduit en tant que k-module.

**4.5.** Avant d'introduire le produit tensoriel réduit pour des espaces, il nous faut le définir dans un contexte simplicial, soit

$$ {}_{S_1}\mathcal{D}^*(\Delta_{s_1})\ \overline{\otimes}\ ...\ \overline{\otimes}\ {}_{S_m}\mathcal{D}^*(\Delta_{s_m}) $$

Pour simplifier les notations, nous nous limiterons à deux facteurs seulement, soit ${}_S\mathcal{D}^*(\Delta_s)\ \overline{\otimes}\ {}_T\mathcal{D}^*(\Delta_t)$, la généralisation à m facteurs ne présentant pas de difficulté particulière.

D'après 2.8, nous pouvons écrire ${}_S\mathcal{D}^*(\Delta_n)$ comme un quotient de ${}_S\mathcal{D}^*(k[x_0, ..., x_n])$, noté simplement ${}_S\mathcal{D}^*(K_{n+1})$, où $K_{n+1}$ est un cube de dimension n+1. Le produit tensoriel réduit précédent est alors l'image du sous-module

$$ {}_S\mathcal{D}^*(K_{s+1})\ \overline{\otimes}\ {}_T\mathcal{D}^*(K_{t+1}) \text{ de } {}_S\mathcal{D}^*(K_{s+1}) \otimes {}_T\mathcal{D}^*(K_{t+1}) $$

par l'homomorphisme de projection

$$ {}_S\mathcal{D}^*(K_{s+1}) \otimes {}_T\mathcal{D}^*(K_{t+1}) \longrightarrow {}_S\mathcal{D}^*(\Delta_s) \otimes {}_T\mathcal{D}^*(\Delta_t) $$

De cette manière, ${}_S\mathcal{D}^*(\Delta_s)\ \overline{\otimes}\ {}_T\mathcal{D}^*(\Delta_t)$ est interprété comme le quotient de ${}_S\mathcal{D}^*(K_{s+1})\ \overline{\otimes}\ {}_T\mathcal{D}^*(K_{t+1})$ par le sous-module engendré par $I_s\ \overline{\otimes}\ {}_S\mathcal{D}^*(K_{t+1})$ et ${}_T\mathcal{D}^*(K_{s+1})\ \overline{\otimes}\ I_t$ avec les notations de 2.8. Ce sous-module est engendré librement par les produits tensoriels de fonctions de Dirac ou de fonctions de Heaviside

$$ (f_0 \otimes ... \otimes f_s) \otimes (g_0 \otimes ... \otimes g_t) $$

les supports singuliers de $f_i$ et $g_j$ étant disjoints pour tout couple (i, j) : on explicite en fait un produit réduit ${}_S\mathcal{D}^*(\Delta_s)\ \overline{x}\ {}_T\mathcal{D}^*(\Delta_t)$ qui engendre le produit tensoriel réduit.

L'opérateur de projection décrit en 4.4 sur ${}_S\mathcal{D}^*(K_{s+1}) \otimes {}_T\mathcal{D}^*(K_{t+1})$ laisse stable ce sous-module. Il en résulte un opérateur de projection de ${}_S\mathcal{D}^*(\Delta_s) \otimes {}_T\mathcal{D}^*(\Delta_t)$ sur le produit tensoriel réduit ${}_S\mathcal{D}^*(\Delta_s)\ \overline{\otimes}\ {}_T\mathcal{D}^*(\Delta_t)$ compatible avec la différentielle cohomologique[10]. En particulier, la cohomologie de ce produit tensoriel réduit est triviale, sauf en dimension 0 où elle est égale à k.

Il convient de noter que les correspondances

[10] Noter cependant que cet opérateur n'est pas compatible avec les opérations faces.

$$(s_1, ..., s_m) \mapsto {}_{S_1}\mathcal{A}^*(\Delta_{s_1}) \otimes ... \otimes {}_{S_m}\mathcal{A}^*(\Delta_{s_m})$$

$$(s_1, ..., s_m) \mapsto {}_{S_1}\mathcal{A}^*(\Delta_{s_1}) \mathbin{\overline{\otimes}} ... \mathbin{\overline{\otimes}} {}_{S_m}\mathcal{A}^*(\Delta_{s_m})$$

définissent des k-modules m-simpliciaux en un sens évident, en raison des considérations strictes imposées sur les supports singuliers. De même, la correspondance

$$(s_1, ..., s_m) \mapsto {}_{S_1}\mathcal{A}^*(\Delta_{s_1}) \mathbin{\bar{\text{x}}} ... \mathbin{\bar{\text{x}}} {}_{S_m}\mathcal{A}^*(\Delta_{s_m})$$

définit un ensemble m-simplicial.

**4.6.** Si S est fini, rappelons que ${}_S\mathcal{A}^*(X)$ est le k-module des applications simpliciales "finies" de $X_\natural$ dans ${}_S\mathcal{A}^*(\Delta_\natural)$. Par ailleurs, $\mathcal{A}^*(X)$ est défini comme la limite inductive des ${}_{S_n}\mathcal{A}^*(X)$ lorsque $S_n = [-n, +n]$ par exemple. De même, si S et T sont finis, ${}_S\mathcal{A}^*(X) \otimes {}_T\mathcal{A}^*(Y)$ s'identifie au k-module des applications **bisimpliciales** "finies" de $X_\natural$ x $Y_\#$ dans ${}_S\mathcal{A}^*(\Delta_\natural) \otimes {}_T\mathcal{A}^*(\Delta_\#)$. Ceci résulte de la correspondance de Dold-Kan appliquée aux k-modules bisimpliciaux et de l'isomorphisme

$$\mathrm{Hom}(E \otimes E', F \otimes F') \cong \mathrm{Hom}(E, F) \otimes \mathrm{Hom}(E', F')$$

avec les conditions de finitude imposées. Nous définissons maintenant le produit tensoriel réduit ${}_S\mathcal{A}^*(X) \mathbin{\overline{\otimes}} {}_T\mathcal{A}^*(Y)$ comme le sous-module de ${}_S\mathcal{A}^*(X) \otimes {}_T\mathcal{A}^*(Y)$ formé des applications finies de $X_\natural$ x $Y_\#$ dans le k-module bisimplicial ${}_S\mathcal{A}^*(\Delta_\natural) \mathbin{\overline{\otimes}} {}_T\mathcal{A}^*(\Delta_\#)$. Si S et T sont infinis, nous définissons ${}_S\mathcal{A}^*(X) \mathbin{\overline{\otimes}} {}_T\mathcal{A}^*(Y)$ comme la limite inductive des ${}_{S'}\mathcal{A}^*(X) \mathbin{\overline{\otimes}} {}_{T'}\mathcal{A}^*(Y)$, S' et T' parcourant les sous-ensembles finis de S et T respectivement. Nous nous proposons de montrer que l'inclusion

$${}_S\mathcal{A}^*(X) \mathbin{\overline{\otimes}} {}_T\mathcal{A}^*(Y) \longrightarrow {}_S\mathcal{A}^*(X) \otimes {}_T\mathcal{A}^*(Y)$$

est un quasi-isomorphisme si S et T sont infinis. Pour cela, quelques lemmes préliminaires sont nécessaires. Ils permettront aussi de démontrer un résultat plus général pour des ensembles simpliciaux X et Y de dimension finie, si S et T ont suffisamment d'éléments.

Nous fixerons les notations suivantes : ${}_S\overline{\mathcal{A}}^*(\Delta_n)$ est l'algèbre des formes différentielles sur $\Delta_n$ qui s'annulent sur toutes les faces, sauf éventuellement la dernière ; comme il est bien connu, ${}_S\overline{\mathcal{A}}^*(\Delta_n)$ s'identifie aussi au quotient de ${}_S\mathcal{A}^*(\Delta_n)$ par le sous k-module engendré par les formes venant de simplexes dégénérés. Nous définisssons ${}_S\overline{\mathcal{A}}^*(\Delta_n) \mathbin{\overline{\otimes}}_T \overline{\mathcal{A}}^*(\Delta_m)$ comme l'image de ${}_S\mathcal{A}^*(\Delta_n) \mathbin{\overline{\otimes}}_T \mathcal{A}^*(\Delta_m)$ dans ${}_S\overline{\mathcal{A}}^*(\Delta_n) \otimes_T \overline{\mathcal{A}}^*(\Delta_m)$ par la projection canonique ; c'est aussi un sous k-module de ${}_S\mathcal{A}^*(\Delta_n) \mathbin{\overline{\otimes}} {}_T\mathcal{A}^*(\Delta_m)$, car la normalisation dans les deux facteurs n'augmente pas les supports singuliers.

.

**4.7. LEMME.** *Soient* S *et* T *deux sous ensembles (non vides) de* **Z** *tels qu'il existe* t *appartenant à* T *et n'appartenant pas à* S. *Soit* $\sigma = \sum_i \omega_i \otimes \theta_i$ *un élément de* ${}_S\overline{\mathcal{D}}^*(\Delta_n) \,\overline{\otimes}\, {}_T\overline{\mathcal{D}}^*(\Delta_m)$, *dont la restriction à* ${}_S\overline{\mathcal{D}}^*(\Delta_n) \,\overline{\otimes}\, {}_T\overline{\mathcal{D}}^*(\Delta_{m-1})$ *est égale à* 0 (par l'opérateur "dernière face"). *Il existe alors un élément*

$$\overline{\sigma} \in {}_S\overline{\mathcal{D}}^*(\Delta_n) \,\overline{\otimes}\, {}_T\overline{\mathcal{D}}^*(\Delta_{m+1})$$

*dont la restriction à la dernière face de* ${}_S\overline{\mathcal{D}}^*(\Delta_n) \,\overline{\otimes}\, {}_T\overline{\mathcal{D}}^*(\Delta_m)$ *est égale à* σ.

*Démonstration.* Soit Y la fonction de Heaviside de support singulier {t} Alors l'expression

$$\overline{\sigma} = \sum_i \omega_i \otimes \theta_i \, Y(x_{m+1})$$

est l'extension requise.

**4.8. LEMME.** *Soient* S *et* T *deux sous-ensembles de* **Z** *vérifiant l'hypothèse du lemme précédent, ainsi que l'hypothèse symétrique* : *il existe* s *appartenant à* S *et n'appartenant pas à* T. *Soit* $\sigma = \sum_i \omega_i \otimes \theta_i$ *un élément de* ${}_S\overline{\mathcal{D}}^*(\Delta_n) \,\overline{\otimes}\, {}_T\overline{\mathcal{D}}^*(\Delta_m)$, *dont la restriction aux dernières "faces" par rapport aux deux variables dans les deux k-modules*

$${}_S\overline{\mathcal{D}}^*(\Delta_n) \,\overline{\otimes}\, {}_T\overline{\mathcal{D}}^*(\Delta_{m-1}) \text{ et } {}_S\overline{\mathcal{D}}^*(\Delta_{n-1}) \,\overline{\otimes}\, {}_T\overline{\mathcal{D}}^*(\Delta_m)$$

*sont égales à* 0 . *Il existe alors un élément* $\overline{\sigma}'$ *dans* ${}_S\overline{\mathcal{D}}^*(\Delta_{n+1}) \,\overline{\otimes}\, {}_T\overline{\mathcal{D}}^*(\Delta_m)$ *dont la restriction à la "face 2*$^e$ *variable"*

$${}_S\overline{\mathcal{D}}^*(\Delta_{n+1}) \,\overline{\otimes}\, {}_T\overline{\mathcal{D}}^*(\Delta_{m-1})$$

*est égale à 0 et dont la restriction à la "face 1*$^{ère}$ *variable"*

$${}_S\overline{\mathcal{D}}^*(\Delta_n) \,\overline{\otimes}\, {}_T\overline{\mathcal{D}}^*(\Delta_m)$$

*est égale à* σ.

*Démonstration.* Comme dans la démonstration précédente, on pose

$$\overline{\sigma}' = \sum_i \omega_i \, Y(x_{m+1}) \otimes \theta_i$$

où $Y(x_{m+1})$ est cette fois-ci une fonction de Heaviside de support singulier {s}.

**4.9.** Si S et T sont finis et si nous suivons le schéma tracé en 3.8, nous pouvons interpréter le produit tensoriel

$${}_S\mathcal{D}^*(X) \otimes {}_T\mathcal{D}^*(Y)$$

comme un "coproduit" dans un contexte bisimplicial. Avec des définitions évidentes, c'est le coproduit suivant

$$[C^{\#}(X) \otimes C^{\natural}(Y)] \,\Delta\, [{}_S\mathcal{D}^*(\Delta_{\#}) \otimes {}_T\mathcal{D}^*(\Delta_{\natural})]$$

$$\cong \mathrm{Hom}_f({}_S\mathcal{D}_*(\Delta^{\#}) \otimes {}_T\mathcal{D}_*(\Delta^{\natural}), C^{\#}(X) \otimes C^{\natural}(Y))$$

.

Le produit tensoriel réduit

$${}_S\mathcal{D}^*(X)\ \overline{\otimes}\ {}_T\mathcal{D}^*(Y)$$

s'identifie de même au coproduit tensoriel suivant

$$[C^{\#}(X) \otimes C^{\natural}(Y)]\ \Delta\ [{}_S\mathcal{D}^*(\Delta_{\#})\ \overline{\otimes}\ {}_T\mathcal{D}^*(\Delta_{\natural})]$$
$$\cong \mathrm{Hom}_f({}_S\mathcal{D}_*(\Delta^{\#})\ \overline{\otimes}\ {}_T\mathcal{D}_*(\Delta^{\natural}), C^{\#}(X) \otimes C^{\natural}(Y))$$

**4.10**. De manière plus générale, considérons deux k-modules cosimpliciaux plats $M^{\#}$ et $N^{\natural}$ de complexes normalisés associés $\overline{M}^{\#}$ et $\overline{N}^{\natural}$. Considérons les foncteurs associés

$$F(M, N) = [M^{\#} \otimes N^{\natural}]\ \Delta\ [{}_S\mathcal{D}^*(\Delta_{\#})\ \overline{\otimes}\ {}_T\mathcal{D}^*(\Delta_{\natural})] \text{ et}$$
$$G(M, N) = [M^{\#} \otimes N^{\natural}]\ \Delta\ [{}_S\mathcal{D}^*(\Delta_{\#}) \otimes {}_T\mathcal{D}^*(\Delta_{\natural})]$$

Par abus d'écriture on posera aussi

$$F(M, N) = {}_S\mathcal{D}^*(M)\ \overline{\otimes}\ {}_T\mathcal{D}^*(N) \text{ et } G(M, N) = {}_S\mathcal{D}^*(M) \otimes {}_T\mathcal{D}^*(N)$$

D'après le théorème de Dold-Kan appliqué aux bicomplexes, nous pouvons interpréter F(M, N) et G(M, N) comme des foncteurs des complexes normalisés $\overline{M}^{\#}$ et $\overline{N}^{\#}$ respectivement. Plus précisément, en conservant les mêmes notations avec les complexes normalisés de cochaînes et de chaînes, les foncteurs précédents s'écrivent

$$F(M, N) = [\overline{M}^{\#} \otimes \overline{N}^{\natural}]\ \Delta\ [{}_S\overline{\mathcal{D}}^*(\Delta_{\#})\ \overline{\otimes}\ {}_T\overline{\mathcal{D}}^*(\Delta_{\natural})] \text{ et}$$
$$G(M, N) = [\overline{M}^{\#} \otimes \overline{N}^{\natural}]\ \Delta\ [{}_S\overline{\mathcal{D}}^*(\Delta_{\#}) \otimes {}_T\overline{\mathcal{D}}^*(\Delta_{\natural})]$$

Les deux foncteurs F et G s'interprètent comme des cycles "doubles" Z'(Z"(M, N)) en bidegré (0, 0) dans les bicomplexes

$$\tilde{F}(M, N) = [\overline{M}^{\#} \otimes \overline{N}^{\natural}] \otimes [{}_S\overline{\mathcal{D}}^*(\Delta_{\#})\ \overline{\otimes}\ {}_T\overline{\mathcal{D}}^*(\Delta_{\natural})] \text{ et}$$
$$\tilde{G}(M, N) = [\overline{M}^{\#} \otimes \overline{N}^{\natural}] \otimes [{}_S\overline{\mathcal{D}}^*(\Delta_{\#}) \otimes {}_T\overline{\mathcal{D}}^*(\Delta_{\natural})]$$

**4.11. DEFINITION.** On dit que S et T sont transversaux s'ils vérifient les conditions du lemme 4.8 : il existe $s \in S$ (resp. $t \in T$) tel que $s \notin T$ (resp. $t \notin S$).

**4.12. THEOREME.** *Le foncteur* $G(M, N) = {}_S\mathcal{D}^*(M) \otimes {}_T\mathcal{D}^*(N)$ *est exact par rapport à* M *et* N. *Si* S *et* T *sont transversaux, il en est de même du foncteur* F(M, N) = ${}_S\mathcal{D}^*(M)\ \overline{\otimes}\ {}_T\mathcal{D}^*(N)$.

*Démonstration.* Nous ferons seulement la seconde vérification qui est la plus délicate en nous inspirant de la démonstration du théorème 3.3. Elle nécessite plusieurs étapes :

a) Compte tenu des conditions de platitude, le foncteur $\tilde{F}(M, N)$ est exact par rapport à M et N.

b) On montre ensuite que $Z''(\overline{M}, \overline{N})$ est un foncteur exact en $\overline{M}$ et $\overline{N}$. Ceci est une conséquence de la nullité du groupe d'homologie par rapport à la seconde variable, soit $H''(\overline{M}, \overline{N})$, d'après le raisonnement fait en 3.3. Pour démontrer cette nullité, on raisonne par récurrence sur la taille des complexes $\overline{M}$ et $\overline{N}$ et on se ramène ainsi au cas où $\overline{M}$ et $\overline{N}$ sont concentrés en degrés n et m respectivement. Pour des raisons de platitude, on peut même supposer que $\overline{M}$ et $\overline{N}$ sont tous deux égaux à k. La vérification à effectuer est alors exactement exprimée par le lemme 4.7. On notera que $Z''(\overline{M}, \overline{N})$ est un k-module plat en général.
c) Démontrons maintenant que $Z'(Z''(M, N))$ est un foncteur exact en M et N. Par le même raisonnement, on doit montrer que $H'(Z''(M, N)) = 0$, vérification qu'il suffit là aussi d'effectuer pour M et N égaux à k en degrés n et m. Le lemme 4.8 répond alors à la question.

**4.13. THEOREME.** *Supposons que* S *et* T *soient transversaux. L'inclusion naturelle*

$${}_S\mathcal{D}^*(M) \mathbin{\overline{\otimes}} {}_T\mathcal{D}^*(N) \longrightarrow {}_S\mathcal{D}^*(M) \boxtimes {}_T\mathcal{D}^*(N)$$

*est alors un quasi-isomorphisme. En particulier, si* X *et* Y *sont deux ensembles simpliciaux quelconques*, l'*inclusion canonique*

$${}_S\mathcal{D}^*(X) \mathbin{\overline{\otimes}} {}_T\mathcal{D}^*(Y) \longrightarrow {}_S\mathcal{D}^*(X) \boxtimes {}_T\mathcal{D}^*(Y)$$

*est un quasi- isomorphisme si* S *et* T *sont transversaux.*
*Démonstration.* Puisque les foncteurs F et G sont exacts, il suffit de démontrer le théorème lorsque M et N sont concentrés en degrés p et q respectivement et égaux à k. Par ailleurs, le théorème est trivial pour p ou q = 0. On est donc ramené à démontrer le théorème lorsque X (resp. Y) est la sphère $\Delta_p/\partial\Delta_p$ (resp. $\Delta_q/\partial\Delta_q$)[11]. Pour toute sphère $\Delta_n/\partial\Delta_n$ on a en général le diagramme cartésien

$$\begin{array}{ccc} C^{\natural}(\Delta_n/\partial\Delta_n) & \longrightarrow & C^{\natural}(\Delta_n) \\ \downarrow & & \downarrow \\ k & \longrightarrow & C^{\natural}(\partial\Delta_n) \end{array}$$

Puisque F et G sont des foncteurs exacts (une autre manière d'exprimer l'axiome de Mayer-Vietoris), on voit qu'il suffit de vérifier le théorème pour X et Y des simplexes standards, comme nous l'avons montré en 4.5.

**4.14.** De manière analogue, on dit que des sous-ensembles $S_1, ..., S_r$ de **Z** sont transversaux si chaque $S_i$ est transversal à la réunion des autres. Par la même méthode, on monre alors que l'inclusion

$${}_{S_1}\mathcal{D}^*(X_1) \mathbin{\overline{\otimes}} \ldots \mathbin{\overline{\otimes}} {}_{S_r}\mathcal{D}^*(X_r) \longrightarrow {}_{S_1}\mathcal{D}^*(X_1) \boxtimes \ldots \boxtimes {}_{S_r}\mathcal{D}^*(X_r)$$

est un quasi-isomorphisme. Puisque $\mathcal{D}^*(X)$ est la limite inductive des ${}_S\mathcal{D}^*(X)$ lorsque S parcourt les sous-ensembles finis de **Z** et que **Z** x **Z** x ... x **Z** est la réunion de produits $S_1$ x $S_2$ x ... x $S_n$ avec des $S_i$ transversaux, il en résulte, par passage à la limite inductive, le quasi-isomorphisme recherché

$$\mathcal{D}^*(X_1) \mathbin{\overline{\otimes}} \ldots \mathbin{\overline{\otimes}} \mathcal{D}^*(X_r) \longrightarrow \mathcal{D}^*(X_1) \boxtimes \ldots \boxtimes \mathcal{D}^*(X_r)$$

[11] Ce qui correspond à un complexe concentré en degrés 0 et n avec des k-modules isomorphes.

**4.15. THEOREME.** *On a un diagramme commutatif*

$$\begin{array}{ccc} \mathcal{D}^*(X_1)\,\overline{\otimes}\,\dots\,\overline{\otimes}\,\mathcal{D}^*(X_r) & \longrightarrow & \mathcal{D}^*(X_1 \times \dots \times X_r) \\ \downarrow & & \downarrow \\ \mathcal{D}^*(X_{\sigma(1)})\,\overline{\otimes}\,\dots\,\overline{\otimes}\,\mathcal{D}^*(X_{\sigma(r)}) & \longrightarrow & \mathcal{D}^*(X_{\sigma(1)} \times \dots \times \mathcal{D}^*(X_{\sigma(r)}) \end{array}$$

*où les flèches verticales sont induites par la permutation* $\sigma$ *des facteurs du produit tensoriel ou des espaces.*

*Démonstration* (pour deux facteurs seulement, la généralisation à r facteurs étant évidente). Il s'agit de montrer la commutativité du diagramme

$$\begin{array}{ccc} \mathcal{D}^*(X)\,\overline{\otimes}\,\mathcal{D}^*(Y) & \longrightarrow & \mathcal{D}^*(X \times Y) \\ R\downarrow & & \sigma\downarrow \\ \mathcal{D}^*(Y)\,\overline{\otimes}\,\mathcal{D}^*(X) & \longrightarrow & \mathcal{D}^*(Y \times X) \end{array}$$

où R est induit par l'échange des facteurs dans le produit tensoriel et où $\sigma$ est induit par la permutation des espaces X et Y. En effet, ceci est une conséquence de la commutativité du diagramme suivant au niveau des simplexes

$$\begin{array}{ccc} \mathcal{D}^*(\Delta_r)\,\overline{\otimes}\,\mathcal{D}^*(\Delta_r)] & \longrightarrow & \mathcal{D}^*(\Delta_r)\,\overline{\otimes}\,\mathcal{D}^*(\Delta_r)] \\ & \searrow \quad \swarrow & \\ & \mathcal{D}^*(\Delta_r) & \end{array}$$

Ici la flèche horizontale est induite par l'échange des facteurs du produit tensoriel réduit, tandis que les flèches obliques définissent la structure produit dans l'algèbre $\mathcal{D}^*(\Delta_r)$. Notons qu'un point essentiel dans ce raisonnement est le fait que les opérations face et dégénérescence sur chaque facteur du produit tensoriel $\mathcal{D}^*(\Delta_p) \otimes \mathcal{D}^*(\Delta_q)$ respectent le produit tensoriel réduit.

**Note** : **Le théorème suivant achève ainsi la vérification des axiomes**

**4.16.** D'après ce qui précède, le produit des formes différentielles

$$\mu : {}_S\mathcal{D}^*(X) \otimes {}_S\mathcal{D}^*(X) \longrightarrow {}_S\mathcal{D}^*(X)$$

est "commutatif" quand on le restreint au produit tensoriel réduit ${}_S\mathcal{D}^*(X)\,\overline{\otimes}\,{}_S\mathcal{D}^*(X)$. En d'autres termes, la symétrie par rapport aux deux variables sur un élément $\rho$ de ${}_S\mathcal{D}^*(X)\,\overline{\otimes}\,{}_S\mathcal{D}^*(X)$ ne change pas la valeur de $\mu(\rho)$. Cependant, cette propriété est plus faible que celle exprimée en 4.15.

**4.17. LEMME.** *Soit* $\sigma = \sum \omega_i \otimes \theta_i$ *un élément de* ${}_S\overline{\mathcal{D}}^*(\Delta_n)\,\overline{\otimes}\,{}_S\overline{\mathcal{D}}^*(\Delta_m)$, *dont la restriction à* ${}_S\overline{\mathcal{D}}^*(\Delta_n)\,\overline{\otimes}\,{}_S\overline{\mathcal{D}}^*(\Delta_{m-1})$ *est égale à* 0 (par l'opérateur "derniière face"). *Si le cardinal de* S *est plus grand que* n, *il existe alors un élément*

$$\overline{\sigma} \in {}_S\overline{\mathcal{D}}^*(\Delta_n)\,\overline{\otimes}\,{}_S\overline{\mathcal{D}}^*(\Delta_{m+1})$$

*dont la restriction à la dernière face de* ${}_S\overline{\mathcal{D}}^*(\Delta_n)\,\overline{\otimes}\,{}_S\overline{\mathcal{D}}^*(\Delta_m)$ *est égale à* $\sigma$.

*Démonstration.* Notons $P : {}_S\mathcal{D}^*(\Delta_n) \longrightarrow {}_S\overline{\mathcal{D}}^*(\Delta_n)$ la projection canonique. Notons aussi $f_x$ une fonction de Dirac ou de Heaviside de support singulier x (pour la fonction 1, on convient que le support singulier est vide ). Alors ${}_S\overline{\mathcal{D}}^*(\Delta_n)$ est un k-module libre de base les $\omega_i = P(f_{x_0} \otimes f_{x_1} \otimes \dots \otimes f_{x_n})$, avec $x_i \in S$ ou $x_i = \varnothing$ et avec $x_i \neq x_{i+1}$ ; son support singulier est $U_i = \{x_0, ..., x_n\}$ qui est de cardinal au plus n. Avec le choix de tels $\omega_i$ , il s'en suit que la restriction de chaque $\theta_i$ à la dernière face doit être égale à 0. Soit maintenant $Y_i(x_{m+1})$ une fonction de Heaviside par rapport à la (m+1)$^e$-variable dont le support singulier est dans S, mais n'appartient pas à $U_i$. L'expression

$$\overline{\sigma} = \sum_i \omega_i \otimes \theta_i \, Y_i(x_{m+1})$$

est alors l'extension requise.

**4.18. LEMME**. *Soient* m *et* n *deux entiers quelconques et S un sous-ensemble de S ayant un cardinal > n + m + 4. Soit* $\sigma$ *un élément de* ${}_S\overline{\mathcal{D}}^*(\Delta_n) \,\overline{\otimes}\, {}_S\overline{\mathcal{D}}^*(\Delta_m)$, *dont la restriction aux dernières " faces" par rapport aux deux variables dans les deux k-modules*

$${}_S\overline{\mathcal{D}}^*(\Delta_n) \,\overline{\otimes}\, {}_S\overline{\mathcal{D}}^*(\Delta_{m-1}) \text{ et } {}_S\overline{\mathcal{D}}^*(\Delta_{n-1}) \,\overline{\otimes}\, {}_S\overline{\mathcal{D}}^*(\Delta_m)$$

*sont égales à* 0 *. Il existe alors un élément* $\overline{\sigma}'$ *dans* ${}_S\overline{\mathcal{D}}^*(\Delta_{n+1}) \,\overline{\otimes}\, {}_S\overline{\mathcal{D}}^*(\Delta_m)$ *dont la restriction à la "face 2$^e$ variable"*

$${}_S\overline{\mathcal{D}}^*(\Delta_{n+1}) \,\overline{\otimes}\, {}_S\overline{\mathcal{D}}^*(\Delta_{m-1})$$

*est égale à 0 et dont la restriction à la "face 1$^{ère}$ variable"*

$${}_S\overline{\mathcal{D}}^*(\Delta_n) \,\overline{\otimes}\, {}_S\overline{\mathcal{D}}^*(\Delta_m)$$

*est égale à* $\sigma$.

*Démonstration.* On utilise de manière essentielle un résultat établi par M. Zisman (voir l'annexe) qui est le suivant : avec les conditions requises sur S, $\sigma$ peut s'écrire $\sum_i \omega_i \otimes \theta_i$ en sorte que $d\omega_i = 0$, $d\theta_i = 0$, $\omega_i$ et $\theta_i$ ayant des supports singuliers disjoints. Puisque $\omega_i$ ne peut être constant, il en résulte que $\omega_i = d\rho_i$, le support singulier de $\rho_i$ étant disjoint de celui de $\theta_i$. On pose alors

$$\overline{\sigma}' = \sum_i \rho_i \otimes \theta_i$$

**4.19. COROLLAIRE** (théorème de stabilité)**.** *Considérons deux k-modules cosimpliciaux plats* $M^\#$ *et* $N^\natural$ *tels que les complexes normalisés associés* $\overline{M}^\#$ *et* $\overline{N}^\natural$ *soient nuls en degrés >* m *et* n *respectivement. Alors, l'inclusion naturelle*

$${}_S\mathcal{D}^*(M) \,\overline{\otimes}\, {}_S\mathcal{D}^*(N) \longrightarrow {}_S\mathcal{D}^*(M) \otimes {}_S\mathcal{D}^*(N)$$

*est un quasi-isomorphisme si* card(S) > m + n + 4. *En particulier, si* X *et* Y *sont deux ensembles simpliciaux de dimension* ≤ m et n respectivement, *les homomorphismes suivants sont des quasi-isomorphismes pour* card(S) > m + n + 4

$${}_S\mathcal{D}^*(X) \,\overline{\otimes}\, {}_S\mathcal{D}^*(Y) \longrightarrow \mathcal{D}^*(X) \,\overline{\otimes}\, \mathcal{D}^*(Y) \longrightarrow \mathcal{D}^*(X) \otimes \mathcal{D}^*(Y)$$

*Démonstration.* C'est la même que celle du théorème 4.13, compte tenu des lemmes 4.17 et 4.18, qui sont les analogues (pour S = T) des lemmes 4.7 et 4.8.

**4.20. Remarque.** Nous démontrerons un résultat très voisin en 5.12 par une méthode sensiblement différente, utilisant la géométrie des polyèdres.

## 5. Produit tensoriel réduit pour des complexes simpliciaux finis.

**5.1.** Nous allons maintenant changer de point de vue et considérer des complexes simpliciaux finis sur lesquels nous pourrons décrire plus concrètement le produit tensoriel réduit. Ceci nous permettra en outre d'affiner les derniers résultats du paragraphe précédent.

Rappelons d'abord qu'un complexe simplicial est donné par un sous-ensemble X de l'ensemble des parties $\mathcal{P}(I)$ d'un ensemble I quelconque (on suppose que X ne contient pas l'ensemble vide). Si $\sigma$ est un élément de X - appelé "cellule"- et si $\tau$ est un sous-ensemble de $\sigma$, c'est aussi un élément de X (appelé "face" de $\sigma$). Si I est un ensemble ordonné, il est classique d'associer à X un ensemble simplicial noté X', défini de la manière suivante : $X_n$ est l'ensemble des suites $i_0, ..., i_n$ d'éléments de I, telles que $i_0 \leq .... \leq i_n$ et $\{i_0, ..., i_n\}$ un élément de X. Si $\tau$ est une face de $\sigma$, on a un morphisme restriction

$$ {}_S\mathcal{D}^*(\sigma) \longrightarrow {}_S\mathcal{D}^*(\tau) $$

qui n'est autre que l'opérateur face ${}_S\mathcal{D}^*(\Delta_n) \longrightarrow {}_S\mathcal{D}^*(\Delta_p)$ défini en 2.5, pour $p \leq n$, grâce à un changement de variables.

Avec ces notations, ${}_S\mathcal{D}^*(X)$ est défini comme la donnée pour chaque cellule $\sigma$ de X d'un élément $\omega_\sigma$ tel que si $\tau$ est une face de $\sigma$, on a $\omega_{\sigma}|_{\tau} = \omega_\tau$ . Autrement dit, ${}_S\mathcal{D}^*(X)$ est l'égalisateur des deux flèches évidentes

$$ \prod_{\sigma} {}_S\mathcal{D}^*(\sigma) \rightrightarrows \prod_{\tau \subset \sigma} {}_S\mathcal{D}^*(\tau) $$

Si Y est un sous-complexe de X, le morphisme de restriction

$$ {}_S\mathcal{D}^*(X) \longrightarrow {}_S\mathcal{D}^*(Y) $$

est induit sur chaque face $\sigma$ par le morphisme

$$ {}_S\mathcal{D}^*(\sigma) \cong {}_S\mathcal{D}^*(\Delta_n) \longrightarrow {}_S\mathcal{D}^*(\Delta_p) \cong {}_S\mathcal{D}^*(\tau) $$

associé à une application strictement croissante de $\{0, ..., p\}$ dans $\{0, ..., n\}$.

**5.2.** Le produit tensoriel ${}_S\mathcal{D}^*(X) \otimes {}_S\mathcal{D}^*(Y)$ (et plus généralement ${}_S\mathcal{D}^*(X) \otimes {}_T\mathcal{D}^*(Y)$) se définit par transport de structure à partir de ${}_S\mathcal{D}^*(X') \otimes {}_S\mathcal{D}^*(Y')$ : un élément de ce k-module s'exprime comme une combinaison linéaire finie

$$ \sum_{r,s} \lambda_{r,s}\, \omega^r \otimes \theta^s $$

où $\omega^r \in {}_S\mathcal{D}^*(X)$ et $\theta^s \in {}_S\mathcal{D}^*(Y)$. Pour tout couple $(\sigma, \tau)$ de cellules dans X et Y, cette expression s'interprète formellement comme une somme finie

$$\sum_{r,s} \lambda_{r,s}\, \omega^r(\sigma) \otimes \theta^s(\tau)$$

Avec ce dictionnaire, le produit tensoriel réduit $_S\mathcal{D}^*(X) \overline{\otimes}\ _S\mathcal{D}^*(Y)$ s'interprète comme le sous k-module formé des combinaisons linéaires précédentes, où $\omega^r(\sigma)$ et $\theta^s(\tau)$ ont des supports singuliers disjoints par rapport à **toutes** les variables $x_i$ figurant dans l'expression de $\omega^r(\sigma)$ et $\theta^s(\tau)$.

**5.3. PROPOSITION.** *Soit* $\omega$ *un élement de* $_S\mathcal{D}^*(\partial\Delta_n)$, *où* $\Delta_n$ *est le* n-*simplexe, défini pour chaque i par la notation* $\omega_i(x_0,..., \hat{x}_i, ...,x_n)$ *qui est sa restriction à la* i-*face. Notons* $\omega_{i,j}$ *la restriction de* $\omega_i$ *à la j-face*, $\omega_{i,j,k}$ *sa restriction à l'intersection des faces numérotées* (i, j, k), *etc. Alors* $\omega$ *est la restriction de l'élément suivant* $\theta$ *de* $_S\mathcal{D}^*(\Delta_n)$, *où* Y *est une fonction de Heaviside de support singulier quelconque*[12]

$$\theta = \sum_{i} Y(x_i)\, \omega_i - \sum_{\{i,j\}} Y(x_i)Y(x_j)\, \omega_{ij} + \sum_{\{i,j,k\}} Y(x_i)Y(x_j)Y(x_k)\, \omega_{ijk} + ...$$

*Démonstration.* Calculons explicitement la restriction de $\theta$ à la i-face. On trouve (pour i fixé dans les symboles de sommation) :

$$\omega_i + \sum_{\{i,j\}} Y(x_j)\, \omega_{ij} - \sum_{\{i,j\}} Y(x_j)\, \omega_{ij} + \sum_{\{i,j,k\}} Y(x_j)Y(x_k)\, \omega_{ijk} - \sum_{\{i,j,k\}} Y(x_j)Y(x_k)\, \omega_{ijk} + ... = \omega_i$$

**5.4.** La formule précédente peut être généralisée de la manière suivante. Une face L de dimension quelconque de $\Delta_n$ est définie en omettant un certain nombre de coordonnées. Par exemple, la forme différentielle $\omega(x_0, \hat{x}_1, \hat{x}_2, ..., x_n)$ est définie sur la (1,2)-face de $\Delta_n$. On écrira alors

$$Y(x_1)Y(x_2)\, \omega(x_0, \hat{x}_1, \hat{x}_2, ..., x_n)$$

pour son extension à $\Delta_n$. Si on pose $Y_L$ pour $Y(x_1)Y(x_2)$ et $\omega_L$ pour $\omega(x_0, \hat{x}_1, \hat{x}_2, ..., x_n)$, l'expression précédente prend la forme $Y_L\, \omega_L$, écriture qui se généralise de manière évidente. La proposition suivante se démontre alors de manière analogue.

**5.5. PROPOSITION.** *Soient* $L_1, ..., L_r$ *des faces de* $\Delta_n$ *et* $\omega$ *une forme différentielle définie sur la réunion de ces faces. Alors* $\omega$ *est la restriction à cette réunion de la forme différentielle suivante* $\theta$ (dite extension "canonique" mais dépendant du choix de Y) *définie sur* $\Delta_n$ *tout entier*

$$\theta = \sum_{i} Y_{L_i}\, \omega_{L_i} - \sum_{\{i,j\}} Y_{L_i} Y_{L_j}\, \omega_{L_i \cap L_j} + \sum_{\{i,j,k\}} Y_{L_i} Y_{L_j} Y_{L_k}\, \omega_{L_i \cap L_j \cap L_k} - ....$$

**5.6. PROPOSITION**. *Soit* K *un complexe simplicial et soit* L *un sous-complexe. Si* $\omega$ *est une forme différentielle sur* L, *il existe alors une forme différentielle* $\theta$ *sur* K *dont la restriction à* L *est égale à* $\omega$. *En outre, le support singulier de* $\theta$ *est égal à celui de* $\omega$ .

---

[12] La notation $\{i_1, ..., i_k\}$ désigne un ensemble de cardinal k formé des éléments **disjoints** $i_1, ..., i_k$ .

*Démonstration.* Elle est évidente si le support singulier de ω est vide. Dans le cas contraire, on commence par remarquer que le complexe K peut être construit à partir de L par recollement de cellules. Il suffit donc de démontrer la proposition lorsque $K = L \cup \Delta_n$. Puisque $L \cap \Delta_n$ est réunion de cellules de $\Delta_n$, on applique la proposition précédente en choisissant une fonction de Heaviside dont le support singulier appartient au support singulier de ω.

Le lemme technique suivant va jouer un rôle important pour le théorème de stabilité.

**5.7. LEMME.** *Soit* K *un complexe simplicial et soit* L *un sous-complexe. On suppose que* K *est réunion de* r *simplexes maximaux* $K_1, K_2 ..., K_r$ *de dimensions* ≤ n. *Tout élément* ω *de* $_S\mathcal{D}^*(K)$, *localement constant sur* L, *s'écrit alors comme une somme* $\sum_i \omega_i$ *de formes différentielles telles que*

1) *Le support singulier des* $\omega_i$ *est inclus dans le support singulier de* ω *(sur chaque cellule de K)*, $\omega_i|_L = 0$ pour $i > 0$ *et* $\omega_0|_L$ *est localement constante.*

2) *Le cardinal du support singulier*[13] *de* $\omega_i$ *est* $\leq r+n-1$.

*Démonstration.* Nous allons procéder par récurrence sur r. Si r = 1, tout élément ω de $\mathcal{D}^*(K) = \mathcal{D}^*(\Delta_p)$ avec $p \leq n$, est soit un scalaire, soit une somme linéairement indépendante de produits $\omega_i$ de fonctions de Heaviside et de fonctions de Dirac (au plus p). Dans le deuxième cas, le support singulier Ss(ω) n'est pas vide et on peut choisir une fonction de Heaviside Y de support singulier un élément de Ss(ω). Soit $\overline{\omega}_i$ l'extension "canonique" à K, décrite dans 5.5 de la restriction de $\omega_i$ à L. Soit enfin $\bar{c}_L$ l'extension canonique à K de $\omega|_L = c_L$.

Alors $\omega = \sum_{i>0} (\omega_i - \overline{\omega}_i) + (\omega_0 - \overline{\omega}_0 + \bar{c}_L)$ est la décomposition requise (car la somme des $\omega_i|_L$ est égale à $c_L$).

Pour passer de r à r+1 dans la récurrence, posons $K = K_1 \cup ... \cup K_{r+1}$ et $L = L_1 \cup ... \cup L_{r+1}$ avec $L_i = K_i \cap L$. Posons aussi $K' = K_1 \cup ... \cup K_r$ et $L' = L_1 \cup ... \cup L_r$. Soit enfin ω un élément de $\mathcal{D}^*(K)$ qu'on écrira comme un couple[14] (σ, τ), où σ (resp. τ) est la restriction de ω à K' (resp. $K_{r+1}$).

Deux cas peuvent se présenter.

**1er cas : la restriction de σ à ($K' \cap K_{r+1}$) est une fonction constante $c_0$**

Remplaçons alors L' par $L'' = L' \cup (K' \cap K_{r+1})$. D'après l'hypothèse de récurrence, σ peut s'écrire $\sum \sigma_i$ , où les $\sigma_i$ vérifient les hypothèses du lemme pour le couple (K', L"). La restriction de $\sigma_i$ à $K' \cap K_{r+1}$ est en particulier égale à 0 pour i > 0. On peut alors écrire

$$\omega = (\sigma, \tau) = \sum_{i>0} (\sigma_i, 0) + (0, \tau - c_0) + (\sigma_0, c_0) \quad (S)$$

[13] On notera Ss(σ) en général le support singulier de σ qui est la **réunion** de tous les supports singuliers nécessités dans chaque cellule.

[14] Plus précisément, la notation (σ, τ) représente un couple constitué d'une forme sur K' et d'une forme sur $K_{r+1}$, coïncidant sur leur intersection.

Par ailleurs, d'après le premier cran de la récurrence, $\tau$ s'écrit aussi comme une somme $\sum \tau_i$ avec $\tau_i \mid L_{r+1} \cup (K' \cap K_{r+1}) = 0$ pour i > 0 et localement constante sur ce sous-complexe pour i = 0.

Deux sous-cas sont alors à examiner :

a) $\tau_0 \mid_{L_{r+1} \cup (K' \cap K_{r+1})}$ est constante (par conséquent égale à $c_0$). Alors chaque terme de la somme précédente a bien un support singulier de cardinal $\leq r + n - 1$ et vérifie les hypothèses du lemme

b) $\theta = \tau_0 \mid_{L_{r+1} \cup (K' \cap K_{r+1})}$ n'est pas constante mais seulement localement constante. Ceci implique qu'il existe un élément s du support singulier de $\tau = \omega|_{K_{r+1}}$ . On peut alors utiliser une fonction de Heaviside ayant cet élément comme support singulier pour étendre "canoniquement" $\theta$ à $K_{r+1}$ tout entier en une fonction $\bar{\theta}$. On réécrit alors la somme (S) précédente sous la forme

$$\omega = (\sigma, \tau) = \sum_{i>0} (\sigma_i, 0) + \sum_{i>0} (0, \tau_i) + (0, \tau_0 - \bar{\theta}) + (\sigma_0, \bar{\theta})$$

Il convient de noter que chaque terme de cette somme a un support singulier de cardinal $\leq$ r+n (c'est le dernier terme qui peut faire augmenter d'un cran le cardinal du support singulier). Par ailleurs, la restriction à L est bien égale à 0, sauf pour le dernier terme qui est localement constant sur L.

**2e cas : la restriction de $\sigma$ à ($K' \cap K_{r+1}$) n'est pas une fonction constante**

Dans ce cas, le support singulier de $\tau = \omega_{\mid K_{r+1}}$ n'est pas vide. Soit Y une fonction de Heaviside dont le support singulier appartient à $Ss(\tau)$.

Pour i > 0, soit $\bar{\sigma}_i$ l'extension canonique à $K_{r+1}$ de $\sigma_i \mid_{K' \cap K_{r+1}}$ et de 0 sur $L_{r+1}$

Soit $\bar{\sigma}_0$ l'extension canonique à $K_{r+1}$ de $\sigma_0 \mid_{K' \cap K_{r+1}}$ et de $\omega$ restreinte à $L_{r+1}$ (qui est localement constante).

On peut alors écrire

$$(\sigma, \tau) = \sum (\sigma_i , \bar{\sigma}_i ) + (0, \tau - \sum \bar{\sigma}_i)$$

Là aussi, le cardinal du support singulier de $(\sigma_i , \bar{\sigma}_i )$ peut augmenter d'un cran par rapport à celui de $\sigma_i$. Les supports singuliers de chaque terme de la somme ont donc bien un cardinal au plus égal à r + n et les hypothèses du lemme sont bien vérifiées.

**5.8. PROPOSITION.** *Soient* $X_1$ *et* $X_2$ *deux complexes simpliciaux finis et soit* $Y_1$ *un sous-complexe de* $X_1$. *On suppose que* $X_2$ *est de dimension ≤ n et peut être recouvert par r simplexes. Alors l'homomorphisme de restriction*

$$ {}_S\mathcal{A}^*(X_1)\,\overline{\otimes}\,{}_S\mathcal{A}^*(X_2) \longrightarrow {}_S\mathcal{A}^*(Y_1)\,\overline{\otimes}\,{}_S\mathcal{A}^*(X_2) $$

*est surjectif si* S *a au moins* r + n *éléments.*

*Démonstration.* Considérons un élément $\sum \omega_i \otimes \theta_i$ du produit tensoriel réduit ${}_S\mathcal{A}^*(Y_1)\,\overline{\otimes}\,{}_S\mathcal{A}^*(X_2)$. D'après le lemme précédent, on ne restreint pas la généralité en supposant que le support singulier de $\theta_i$ contient au plus r + n - 1 éléments. Pour chaque i, il existe donc un élément $s_i$ de S qui n'est pas dans le support singulier de $\theta_i$. La fonction de Heaviside associée permet alors d'étendre canoniquement $\omega_i$ en une forme $\overline{\omega}_i$ sur X dont le support singulier est disjoint de celui de $\theta_i$. L'extension requise est alors $\sum \overline{\omega}_i \otimes \theta_i$ .

La proposition suivante se démontre de manière analogue.

**5.9. PROPOSITION.** *Plus généralement, soient* $X_1, \dots , X_m$ *des complexes simpliciaux de dimension* ≤ n *et pouvant être recouverts par* r *simplexes. Soient* $Y_1, \dots , Y_m$ *des sous-complexes de* $X_1, \dots, X_m$ *respectivement. Alors l'homomorphisme de restriction*

$$ {}_S\mathcal{A}^*(X_1)\,\overline{\otimes}\dots\overline{\otimes}\,{}_S\mathcal{A}^*(X_m) \longrightarrow {}_S\mathcal{A}^*(Y_1)\,\overline{\otimes}\dots\overline{\otimes}\,{}_S\mathcal{A}^*(X_m) $$

*est surjectif si* S *a au moins* (m - 1)(r + n - 1) + 1 *éléments.*

**5.10. THEOREME.** *Sous les hypothèses de la proposition* 5.9, *l'homomorphisme d'inclusion*

$$ {}_S\mathcal{A}^*(X_1)\,\overline{\otimes}\dots\overline{\otimes}\,{}_S\mathcal{A}^*(X_m) \longrightarrow \mathcal{A}^*(X_1) \otimes\dots\otimes {}_S\mathcal{A}^*(X_m) $$

*est un quasi-isomorphisme si* S *a au moins* (m - 1)(r + n - 1) + 1 *éléments*

*Démonstration.* La proposition précédente et un argument de Mayer-Vietoris standard permet de nous ramener au cas où les $X_i$ sont des simplexes et les $Y_i$ des sous-simplexes . L'opérateur de projection

$$ {}_S\mathcal{A}^*(X_1) \otimes\dots\otimes {}_S\mathcal{A}^*(X_m) \longrightarrow \mathcal{A}^*(X_1)\,\overline{\otimes}\dots\overline{\otimes}\,{}_S\mathcal{A}^*(X_m) $$

introduit en 4.3 est compatible avec la différentielle **cohomologique** ; il est l'inverse à gauche de l'inclusion canonique. Par conséquent, ${}_S\mathcal{A}^*(X_1)\,\overline{\otimes}\dots\overline{\otimes}\,{}_S\mathcal{A}^*(X_m)$ a la bonne cohomologie, c'est-à-dire réduite à k en degré 0.

**5.11. THEOREME.** *Soient* M *et* N *deux* k-*modules cosimpliciaux de dimension finie avec* N *de dimension* ≤ n, *dans le sens précisé au § 4, On se place dans l'une des trois hypothèses suivantes*

*a)* k *est un corps*

*b)* k *est un anneau principal*

*c)* M *et* N *sont les modules cosimpliciaux associés aux modules des cochaînes sur des ensembles simpliciaux X et Y respectivement sur un anneau cohérent k quelconque.*

*En supposons vraie l'une de ces trois conditions, l'inclusion canonique*

$$ {}_S\mathcal{A}^*(M)\,\overline{\otimes}\,{}_S\mathcal{A}^*(N) \longrightarrow {}_S\mathcal{A}^*(M) \otimes {}_S\mathcal{A}^*(N) $$

*est alors un quasi-isomorphisme si* S *a au moins* 2n+1 *éléments. En particulier, si* X *et* Y *sont deux ensembles simpliciaux de dimension* finie $\leq$ n, *et avec les mêmes conditions sur* S (k *étant un anneau cohérent quelconque)*, *l'application*

$$ {}_S\mathcal{A}^*(X)\ \overline{\otimes}\ {}_S\mathcal{A}^*(Y) \longrightarrow {}_S\mathcal{A}^*(X) \otimes {}_S\mathcal{A}^*(Y) $$

*est un quasi-isomorphisme.*

*Démonstration.* Supposons d'abord que k soit un corps. Dans ce cas, tout module cosimplicial est somme directe de modules “élémentaires” c'est-à-dire dont les complexes de cochaînes associés (par la correspondance de Dold-Kan) sont soit concentrés en un seul degré, soit du type....$0 \longrightarrow A \xrightarrow{\alpha} B \longrightarrow 0$ ... , où $\alpha$ est un isomorphisme. Le théorème étant évident pour un complexe en degré 0, on est finalement amené à vérifier le quasi-isomorphisme dans un cadre topologique, soit pour l'inclusion

$$ {}_S\mathcal{A}^*(X)\ \overline{\otimes}\ {}_S\mathcal{A}^*(Y) \longrightarrow {}_S\mathcal{A}^*(X) \otimes {}_S\mathcal{A}^*(Y) $$

avec X ou Y = $\Delta_p/\partial\Delta_p$ ou $\Delta_p/\Lambda\Delta_p$,, $\Lambda\Delta_p$ étant le bord de $\Delta_p$ privé d'une face ($p \leq n$)[15]. Vérifions le quasi-isomorphisme par exemple si $X = \Delta_q/\partial\Delta_q$ et $Y = \Delta_p/\partial\Delta_p$ (le raisonnement est analogue dans les autres cas). On a un diagramme cartésien

$$ \begin{array}{ccc} {}_S\mathcal{A}^*(\Delta_m/\partial\Delta_m) & \longrightarrow & {}_S\mathcal{A}^*(\Delta_m) \\ \downarrow & & \downarrow \\ {}_S\mathcal{A}^*(P) & \longrightarrow & {}_S\mathcal{A}^*(\partial\Delta_m) \end{array} $$

où P est un point. Ceci permet d'identifier ${}_S\mathcal{A}^*(\Delta_m/\partial\Delta_m)$ au k-module des formes différentielles sur $\Delta_m$ qui sont constantes sur le bord. Pour démontrer le quasi-isomorphisme recherché

$$ {}_S\mathcal{A}^*(X)\ \overline{\otimes}\ {}_S\mathcal{A}^*(Y) \longrightarrow {}_S\mathcal{A}^*(X) \otimes {}_S\mathcal{A}^*(Y) $$

on écrit la suite exacte

$$ 0 \longrightarrow {}_S\mathcal{A}^*(\Delta_q/\partial\Delta_q)\ \overline{\otimes}\ {}_S\mathcal{A}^*(\Delta_p/\partial\Delta_p) \longrightarrow {}_S\mathcal{A}^*(\Delta_q/\partial\Delta_q)\ \overline{\otimes}\ {}_S\mathcal{A}^*(\Delta_p) $$

$$ \longrightarrow {}_S\mathcal{A}^*(\Delta_q/\partial\Delta_q)\ \overline{\otimes}\ {}_S\mathcal{A}^*(\partial\Delta_p) \longrightarrow 0 $$

L'exactitude de cette suite est évidente, sauf pour la surjectivité à la fin : celle-ci se démontre si card(S) $\geq$ 2q + 1 en suivant le même principe que dans la proposition 5.8. En effet, tout élément de ${}_S\mathcal{A}^*(\Delta_q/\partial\Delta_q)\ \overline{\otimes}\ {}_S\mathcal{A}^*(\partial\Delta_p)$ s'écrit comme une somme $\sum \omega_i \otimes \theta_i$, où les $\theta_i$ sont des éléments de ${}_S\mathcal{A}^*(\Delta_q)$ qui sont constants sur $\partial\Delta_q$ et tels que leurs supports singuliers aient un cardinal $\leq q + 1$.

On démontre de même l'exactitude de la suite

$$ 0 \longrightarrow {}_S\mathcal{A}^*(\partial\Delta_q)\ \overline{\otimes}\ {}_S\mathcal{A}^*(\Delta_p/\partial\Delta_p) \longrightarrow {}_S\mathcal{A}^*(\partial\Delta_q)\ \overline{\otimes}\ {}_S\mathcal{A}^*(\Delta_p) $$

$$ \longrightarrow {}_S\mathcal{A}^*(\partial\Delta_q)\ \overline{\otimes}\ {}_S\mathcal{A}^*(\partial\Delta_p) \longrightarrow 0 $$

si card(S) $\geq$ 2q+1 (car $\partial\Delta_q$ peut être recouvert par q + 1 simplexes). Par ces réductions successives, on voit finalement qu'on se ramène à $X = \partial\Delta_q$ et $Y = \partial\Delta_p$ qui sont des complexes simpliciaux de type q+1 et p+1 respectivement. Le théorème 5.10 permet alors de conclure la démonstration du théorème si k est un corps.

[15] Ce qui correspond aux deux cas envisagés par la correspondance de Dold-Kan de nouveau.

Supposons maintenant que k soit un anneau principal. Pour démontrer le quasi-isomorphisme

$$_S\mathcal{D}^*(M) \mathbin{\overline{\otimes}} {}_S\mathcal{D}^*(N) \longrightarrow {}_S\mathcal{D}^*(M) \otimes {}_S\mathcal{D}^*(N)$$

il suffit de vérifier que le cône C de cette application est acyclique. Puisque C est plat, on a une suite exacte

$$0 \longrightarrow C \longrightarrow C \longrightarrow C/p \longrightarrow 0$$

où (p) est un idéal premier quelconque de C, la première flèche étant la multiplication par p. La cohomologie de C/p étant nulle pour tout p, la cohomologie de C est un module sur le corps des fractions F de k. Puisque $F \otimes_k F \cong F$, on en déduit que la cohomologie de C est réduite à 0.

Considérons maintenant le 3e cas où M = C*(X) et N = C*(Y), X et Y étant des espaces de dimension ≤ n. Les deux modules considérés

$$_S\mathcal{D}^*(M) \mathbin{\overline{\otimes}} {}_S\mathcal{D}^*(N) \text{ et } {}_S\mathcal{D}^*(M) \otimes {}_S\mathcal{D}^*(N)$$

sont alors obtenus par extension des scalaires à partir de k = **Z**. Si on désigne par $C_{\mathbf{Z}}$ le cône corrrespondant à l'inclusion canonique, on a donc $C_k \cong C_{\mathbf{Z}} \otimes k$ avec des notations évidentes. Puisque les cocycles Z, donc les cobords B du complexe $C_{\mathbf{Z}}$, sont des **Z**-modules libres, la tensorisation par k des suites du type

$$0 \longrightarrow Z \longrightarrow C \longrightarrow B \longrightarrow 0$$

permet de montrer que la cohomologie de $C_k$ est nulle également.

**5.12. Généralisation.** Si $M_1$, ..., $M_p$ sont de même des k-modules cosimpliciaux de dimension n, on démontre sous les même hypothèses le quasi-isomorphisme

$$_S\mathcal{D}^*(M_1) \mathbin{\overline{\otimes}} \dots \mathbin{\overline{\otimes}} {}_S\mathcal{D}^*(M_p) \longrightarrow {}_S\mathcal{D}^*(M_1) \otimes \dots \otimes {}_S\mathcal{D}^*(M_p)$$

si S a au moins 2(p-1)n + 1 éléments

## 6. Description des cup i-produits et des opérations de Steenrod dans le cadre des algèbres différentielles graduées quasi-commutatives.

**6.1.** Soit A une ADGQ sur un anneau cohérent k telle que, pour tout n, $A^{\bar{\otimes} n}$ et $A^{\otimes n}/A^{\bar{\otimes} n}$ (donc $A^{\otimes n}$) soient projectifs. Ceci sera le cas si par exemple k est un corps ou si $A = \mathcal{D}^*(X)$. Soit $C(\mathfrak{S}_n)$ la bar résolution du groupe symétrique $\mathfrak{S}_n$ (vue comme un complexe cohomologique en degrés ≤ 0). Nous comptons définir un morphisme[16]

$$\mu : C(\mathfrak{S}_n) \otimes_{\mathfrak{S}_n} A^{\otimes n} \longrightarrow A$$

en sorte que le diagramme suivant commute

$$\begin{array}{ccc} C(\mathfrak{S}_n) \otimes_{\mathfrak{S}_n} A^{\otimes n} & \xrightarrow{\ \mu\ } & A \\ \cup\!\!\uparrow & & \| \\ C(\mathfrak{S}_n) \otimes_{\mathfrak{S}_n} A^{\bar{\otimes} n} & \xrightarrow{\ \bar{\mu}\ } & A \end{array}$$

Ici $\bar{\mu}$ désigne l'homomorphisme composé $C(\mathfrak{S}_n) \otimes_{\mathfrak{S}_n} A^{\bar{\otimes} n} \longrightarrow k \otimes_{\mathfrak{S}_n} A^{\bar{\otimes} n} \longrightarrow A^{\bar{\otimes} n}$, où la première flèche est induite par l'augmentation de $C(\mathfrak{S}_n)$ sur k et la seconde par le cup-produit restreint à $A^{\bar{\otimes} n}$ (qui est équivariant). Nous allons exiger aussi que la composante de $\mu$ de degré 0 par rapport à la première variable, soit

$$\mu_0 : C^0(\mathfrak{S}_n) \otimes_{\mathfrak{S}_n} A^{\otimes n} \cong A^{\otimes n} \longrightarrow A$$

soit le cup-produit usuel. Nous appellerons un tel relevé $\mu$ un “cup-produit supérieur”.

Pour construire le morphisme $\mu$, il est commode de considérer la somme amalgamée $T_A$ définie par le diagramme suivant

$$\begin{array}{ccc} C^0(\mathfrak{S}_n) \otimes_{\mathfrak{S}_n} A^{\bar{\otimes} n} & \longrightarrow & C^0(\mathfrak{S}_n) \otimes_{\mathfrak{S}_n} A^{\otimes n} \\ \downarrow & & \downarrow \\ C(\mathfrak{S}_n) \otimes_{\mathfrak{S}_n} A^{\bar{\otimes} n} & \longrightarrow & T_A \end{array}$$

On a un homomorphisme injectif de $T_A$ dans $E_A = C(\mathfrak{S}_n) \otimes_{\mathfrak{S}_n} A^{\otimes n}$ qui est un quasi-isomorphisme. Considérons alors la suite exacte suivante de complexes en homomorphismes

$$0 \longrightarrow \mathcal{HOM}(E_A/T_A, A) \longrightarrow \mathcal{HOM}(E_A, A) \xrightarrow{\ \sigma\ } \mathcal{HOM}(T_A, A) \longrightarrow 0$$

La surjectivité à droite résulte du fait que $T_A$ est projectif. Par ailleurs, la cohomologie du complexe $\mathcal{HOM}(E_A/T_A, A)$ est triviale car $E_A/T_A$ est acyclique et projectif : il possède donc un opérateur d'homotopie qui permet de démontrer l'acyclicité du complexe $\mathcal{HOM}(E_A/T_A, A)$. Un raisonnement identique à celui fait en 3.3 permet de montrer l'exactitude de la suite

[16] Le produit tensoriel désigne ici le complexe total associé au bicomplexe formé avec les différentielles des deux facteurs.

$$0 \longrightarrow Z^0\mathcal{HOM}(E_A/T_A, A) \longrightarrow Z^0\mathcal{HOM}(E_A, A) \longrightarrow Z^0\mathcal{HOM}(T_A, A) \longrightarrow 0$$

Par conséquent, puisque $\bar{\mu} \in Z^0\mathcal{HOM}(T_A, A)$, il existe bien $\mu \in Z^0\mathcal{HOM}(E_A, A)$ tel que $\sigma(\mu) = \bar{\mu}$.

**6.2.** Une autre manière de voir les choses est de considérer le complexe $\mathcal{HOM}_\natural(A^{\otimes n}, A)$. En degrés $r > 0$, $\mathcal{HOM}_r(A^{\otimes n}, A)$ désigne le k-module des applications k-linéaires de degré -r de $A^{\otimes n}$ dans A, tandis qu'en degré 0, $\mathcal{HOM}_0(A^{\otimes n}, A)$ est le k-module des morphismes entre les complexes $A^{\otimes n}$ et A. Nous avons ainsi un complexe de nature "homologique"

$$\xrightarrow{D} \mathcal{HOM}_2(A^{\otimes n}, A) \xrightarrow{D} \mathcal{HOM}_1(A^{\otimes n}, A)$$
$$\xrightarrow{D} \mathcal{HOM}_0(A^{\otimes n}, A) \longrightarrow 0$$

où D est défini par la formule usuelle

$$D(f)(\omega) = df(\omega) + (-1)^{\deg f} f(d\omega)$$

Il convient de noter que le noyau de D est formé de morphismes de complexes (avec un décalage de degré approprié). Le complexe $\mathcal{HOM}_\natural(A^{\bar{\otimes} n}, A)$ est défini de manière analogue et, d'après ce qui précède, on a un morphisme surjectif

$$\mathcal{HOM}_\natural(A^{\otimes n}, A) \longrightarrow \mathcal{HOM}_\natural(A^{\bar{\otimes} n}, A)$$

ainsi qu'un diagramme commutatif de **complexes** (avec des notations analogues)

$$\begin{array}{ccc} \mathcal{HOM}_\natural(A^{\otimes n}, A) & \longrightarrow & \mathcal{HOM}_\natural(A^{\bar{\otimes} n}, A) \\ & \mu \nwarrow & \uparrow \bar{\mu} \\ & & C_\natural(\mathfrak{S}_n) \end{array}$$

Le morphisme de complexes $\mu = (\mu_n)$ peut alors se construire alors par récurrence sur n en partant du cup-produit usuel $\mu_0$ et en exploitant l'acyclicité du complexe quotient $A^{\otimes n}/A^{\bar{\otimes} n}$.

**6.3.** On peut remplacer le groupe symétrique $\mathfrak{S}_n$ par le groupe cyclique $C_n$ et la bar-résolution $C_\natural(\mathfrak{S}_n)$ par la résolution standard du groupe cyclique, soit

$$M_i \xrightarrow{\alpha_i} M_{i-1} \longrightarrow \ldots M_1 \longrightarrow M_0 \longrightarrow k$$

où $M_i = k[x]/x^n - 1$ et où $\alpha_i$ est la multiplication par 1 - x si i est impair et la multiplication par $1 + x + \ldots + x^{n-1}$ si i est pair. L'image du générateur $x_i = x$ de $E_i$ définit une application k-linéaire de degré -i de $A^{\otimes n}$ vers A, appelée cup i-produit et notée $\mu_i$ qui généralise le cup-produit ordinaire $\mu_0$. Ces homomorphismes $\mu_i$ vérifient les relations suivantes (D étant défini

ci-dessus)

$D\mu_1 = (1 - t)\,\mu_0$

$D\mu_2 = (1 + t + ... + t^{n-1})\,\mu_1$

$D\mu_3 = (1 - t)\,\mu_2$

etc. ,

où t désigne l'action du générateur du groupe cyclique $C_n$ sur le complexe $\mathcal{HOM}_{k}(A^{\otimes n}, A)$. Supposons maintenant que k soit de caractéristique n. On remarque alors que l'image de l'application "puissance $n^{\text{ième}}$"

$$A \longrightarrow A^{\otimes n}$$

est annulée[17] à la fois par 1 - t et par $1 + t + ... + t^{n-1}$. L'application composée

$$A \longrightarrow A^{\otimes n} \xrightarrow{\mu_i} A$$

envoit donc un cocycle en un cocycle, soit de $Z^r(A)$ dans $Z^{nr - i}(A)$. Si n est un nombre premier p, on verra un peu plus loin que pour $A = \mathcal{D}^*(X)$, cette application induit en cohomologie l'opération de Steenrod à une normalisation près

$$D_i : H^r(X) \longrightarrow H^{pr - i}(X)$$

(si p = 2 par exemple on trouve le carré de Steenrod $Sq^j$ avec j = r - i).

**6.4.** Ce cup produit supérieur $\mu_i$ est"naturel" dans un sens que nous allons préciser. Considérons un morphisme surjectif d'ADGQ, soit $A \longrightarrow B$. Par ailleurs, supposons donné un cup-produit supérieur pour B, soit

$$\begin{array}{ccc} C(\mathfrak{S}_n) \otimes_{\mathfrak{S}_n} B^{\otimes n} & \longrightarrow & B \\ \uparrow & & \| \\ T_B & \longrightarrow & B \end{array}$$

Nous allons montrer qu'il existe un cup-produit supérieur pour A tel que les diagrammes suivants commutent

$$\begin{array}{ccc} C(\mathfrak{S}_n) \otimes_{\mathfrak{S}_n} A^{\otimes n} & \longrightarrow & A \\ \downarrow & & \downarrow \\ C(\mathfrak{S}_n) \otimes_{\mathfrak{S}_n} B^{\otimes n} & \longrightarrow & B \end{array} \qquad \begin{array}{ccc} C(\mathfrak{S}_n) \otimes_{\mathfrak{S}_n} A^{\otimes n} & \longrightarrow & A \\ \uparrow & & \| \\ T_A & \longrightarrow & A \end{array}$$

Pour alléger les notations, posons $E = C(\mathfrak{S}_n) \otimes_{\mathfrak{S}_n} A^{\otimes n}$ et $E' = T_A$ qui lui est quasi-isomorphe par l'inclusion de E' dans E. Le problème (classique) revient à construire une flèche de E dans A tel que le diagramme suivant commute

$$\begin{array}{ccc} E' & \longrightarrow & A \\ \downarrow & \nearrow & \downarrow \\ E & \longrightarrow & B \end{array}$$

[17] L'action du groupe cyclique tient compte du degré des cochaînes.

En effet, la commutativité des deux triangles est une reformulation de la commutativité des deux diagrammes ci-dessus.
Si on désigne par U le noyau de la flèche A ⟶ B et par V le conoyau de la flèche E' ⟶ E, on a une suite exacte

$$0 \longrightarrow \mathcal{HOM}(V, U) \rightarrow \mathcal{HOM}(E, A) \rightarrow \mathcal{HOM}(E, B) \times \mathcal{HOM}(E', A) \longrightarrow 0$$

Puisque V est acyclique et projectif, le complexe $\mathcal{HOM}(V, U)$ est acyclique, ce qui démontre l'assertion par un raisonnement analogue à celui fait en 3.3.

**6.5.** Cette naturalité du cup-produit supérieur implique dans le cas topologique une propriété homotopique. En effet, soient

$$\mu \text{ et } \mu' : C(\mathfrak{S}_n) \otimes_{\mathfrak{S}_n} B'^{\otimes n} \longrightarrow B'$$

deux cup-produits supérieurs avec $B' = \mathcal{D}^*(X)$. Si on pose $A = \mathcal{D}^*(X \times [0, 1])$ avec la flèche restriction évidente de A dans $B = B' \times B' = \mathcal{D}^*(X \times \{0\}) \times \mathcal{D}^*(X \times \{1\})$ , on voit que μ et μ' sont homotopes en un sens évident. En particulier l'opération de $Z^r(X)$ dans $Z^{pr-i}(X)$ est invariante par homotopie, donc induit une transformation cohomologique de $H^r(X)$ dans $H^{pr-i}(X)$ qu'on identifiera plus loin à l'opération de Steenrod $D_i$ telle qu'elle est définie dans [10] par exemple.

**6.6.** En fait, de multiples définitions des opérations de Steenrod existent dans la littérature (cf. la bibliographie attachée à [10] par exemple). Elles reposent pour la plupart sur une théorie des cochaînes particulière. Il convient de remarquer cependant que le résultat final au niveau de la cohomologie n'en dépend pas.
Ce dernier point est basé sur le principe suivant. Considérons deux ADG simpliciales A* et B* vérifiant les axiomes de Cartan (cf. [3], ainsi que 2.7/8). Il en est alors de même pour l'ADG simpliciale $A^* \otimes B^*$ , avec des hypothèses de platitude par exemple. Cette situation se présente dans le cas où B* est le complexe des formes différentielles non commutatives (ou des cochaînes normalisées) sur le simplexe standard (cf. [10]). Nous avons ainsi des quasi-isomorphismes entre les théories de cochaînes associées, associés aux morphismes de A* vers $A^* \otimes B^*$ et de B* vers $A^* \otimes B^*$ , définis respectivement par $a \mapsto a \otimes 1$ et $b \mapsto 1 \otimes b$. Ainsi, toutes les opérations cohomologiques définies via A* ou B* peuvent se comparer via $A^* \otimes B^*$.

**6.7.** Supposons donnée ainsi une ADG simpliciale, notée $A_\#^*$ , vérifiant les axiomes standard rappelés en 2.7/8 (les symbole * et # désignant respectivement le degré et la dimension simpliciale). Si X est un ensemble simplicial, A*(X) est l'ensemble des morphismes[18] de X

[18] Pour simplifier, nous pouvons aussi supposer que les morphismes sont "finis" dans le sens que l'image d'un n simplexe est une combinaison linéaire d'élements dégénérés si n est assez grand.

dans A*.

Soit $Z^r$ le k-module simplicial $s \mapsto Z^r(\Delta_s)$, le k-module des cocycles de degré r sur $\Delta_s$. Comme il est montré dans [3], c'est un modèle de l'espace d'Eilenberg-Mac Lane K(k, r). Plus généralement, si X est un ensemble simplicial, l'ensemble simplicial

$$s \mapsto Z^r(X \times \Delta_s)$$

représente l'espace fonctionnel des applications simpliciales de X dans $Z^r$.

**6.8.** L'application "élévation à la puissance n" définit un morphisme de $Z^r$ dans $Z^{r*n}$, où nous désignons par $Z^{r*n}$ l'espace des cycles de degré r.n dans $(A^*)^{\otimes n}$. Cette application est équivariante pour l'action du groupe symétrique $\mathfrak{S}_n$, triviale sur $Z^r$ et permutant les facteurs dans $Z^{r*n} \subset (A^*)^{\otimes n}$. Analysons cette dernière action sur $Z^{r*n}$. D'après la correspondance de Dold-Kan déjà utilisée dans les paragraphes précédent , $Z^{r*n}$ correspond à un certain complexe de chaînes

$$C_{m+1} \longrightarrow C_m \longrightarrow C_{m-1}$$

muni d'une action du groupe symétrique $\mathfrak{S}_n$. On peut alors écrire un diagramme commutatif

$$\begin{array}{ccccccc}
\longrightarrow C_{r.n+1} & \longrightarrow & C_{r.n} & \longrightarrow & C_{r.n-1} & \longrightarrow \\
\| & & \uparrow & & \uparrow & \\
\longrightarrow C_{r.n+1} & \longrightarrow & Z_{r.n} & \longrightarrow & 0 & \\
\downarrow & & \downarrow & & \downarrow & \\
0 & \longrightarrow & H_{r.n} & \longrightarrow & 0 &
\end{array}$$

où $H_{r.n}$ est le k-module k concentré en degré r.n. La traduction simpliciale de ce diagramme commutatif de complexes est donnée par des équivalences d'homotopie équivariantes entre espaces d'Eilenberg-Mac Lane

$$Z^{r*n} \longleftarrow \Gamma_{n,r} \longrightarrow \tilde{Z}^{r.n}$$

où $\tilde{Z}^{r.n}$ est l'espace "minimal", correspondant à un complexe de chaînes concentré en un seul degré. Le groupe symétrique opère sur $\tilde{Z}^{r.n}$ par le signe $(-1)^{r.n}$ pour une transposition. En composant avec la puissance $n^{\text{ième}}$, on a ainsi des morphismes équivariants de k-modules simpliciaux (les deux derniers étant des équivalences d'homotopie) :

$$Z^r \longrightarrow Z^{r.n} \longleftarrow \Gamma_{n,r} \longrightarrow \tilde{Z}^{r.n} \qquad \text{(E)}$$

**6.9.** Deux cas particuliers importants de la suite (E) sont à considérer : $A^* = \Omega^*$ est la théorie des formes différentielles non commutatives étudiée dans [10] ou bien A* est la théorie des formes différentielles quasi-commutatives $\mathfrak{Q}^*$ introduite dans cet article. Dans le premier cas, la suite (E) peut se réduire en une suite à trois termes. En effet, le complexe des formes différentielles antisymétriques (noté $\tilde{\Omega}^n$ car il correspond aux espaces d'Eilenberg-Mac Lane minimaux ; cf. [9]) est quasi-isomorphe à $\Omega^n$ par un morphisme d'inclusion $\tilde{\Omega}^n \rightarrowtail \Omega^n$

On peut donc écrire un diagramme commutatif d'espaces et de morphismes équivariants

$$\begin{array}{ccccccc} Z^r & \longrightarrow & Z^{r*n} & \longleftarrow & \Gamma_{n,r} & \longrightarrow & \tilde{Z}^{r.n} \\ \| & & \| & & \downarrow & & \| \\ Z^r & \longrightarrow & Z^{r.n} & \longleftarrow & \tilde{Z}^{r.n} & = & \tilde{Z}^{r.n} \end{array}$$

On procède de manière analogue pour la théorie $\mathcal{D}^*$, en remplaçant $\tilde{\Omega}^{*n}$ par $\mathcal{D}^{*\overline{\otimes}n}$. Nous obtenons alors le diagramme commutatif suivant :

$$\begin{array}{ccccccc} Z^r & \longrightarrow & Z^{r.n} & \longleftarrow & \Gamma_{n,r} & \longrightarrow & \tilde{Z}^{r.n} \\ & & \uparrow & & \uparrow & & \| \\ & & Z^{r(\overline{\otimes}n)} & \leftarrow & \Gamma'_{n,r} & \longrightarrow & \tilde{Z}^{r.n} \end{array}$$

Dans ce diagramme, on note $Z^{r(\overline{\otimes}n)}$ le k-module des cocycles du complexe $(\mathcal{D}^*)^{\overline{\otimes}n}$ en degré r.n. Il convient de noter que le <u>produit</u> $Z^{r(\overline{\otimes}n)} \longrightarrow Z^{r.n}$ réalise directement une équivalence d'homotopie équivariante entre $Z^{r(\overline{\otimes}n)}$ dans $\tilde{Z}^{r.n}$ car ce produit est commutatif (au signe près). On a ainsi des équivalences d'homotopie en zigzag

$$Z^r \longrightarrow Z^{r.n} \longleftarrow Z^{r(\overline{\otimes}n)} \cong \tilde{Z}^{r.n}$$

**6.10.** Pour simplifier, restreignons nous à un sous-groupe G de $\mathfrak{S}_n$ tel que G opère trivialement sur $\tilde{Z}^{r.n}$ (par exemple le groupe cyclique $C_n$ si n = p est un nombre premier avec un choix de k = $F_p$ , le corps fini à p éléments). Nous noterons EG l'ensemble simplicial acyclique associé au fibre universel EG sur BG. Par abus d'écriture, les k-modules associées seront aussi notés EG et BG.
De la suite vue en 6.8, nous déduisons des morphismes de k-modules simpliciaux bien définis à homotopie près :

$$Z^r \longrightarrow \mathrm{Mor}_G(EG, Z^{r*n}) \longleftarrow \mathrm{Mor}_G(EG, \Gamma_{m,r}) \longrightarrow \mathrm{Mor}_G(EG, \tilde{Z}^{r.n}) = \mathrm{Mor}(BG, \tilde{Z}^{r.n})$$

Cette méthode générale permet de définir une transformation naturelle :

$$H^r(X) \longrightarrow H^{rm}(X \times BG)$$

D'après la remarque faite en 6.6 et la théorie axiomatique développée dans [10], on retrouve bien les opérations de Steenrod classiques, indépendamment de la théorie cohomologique choisie qui vérifie les axiomes de Cartan.

**6.11.** Voici une autre interprétation de cette construction générale faisant spécifiquement appel à la théorie $\mathcal{D}^*$. A partir de maintenant, on notera simplement A* cette théorie. Nous avons en général une équivalence d'homotopie

$$EG \times_G \tilde{Z}^{r.n} \longrightarrow EG \times_G Z^{r.n}$$

et donc une équivalence d'homotopie inverse suivie d'une projection

$$EG \times_G (Z^r)^{\otimes n} \longrightarrow EG \times_G \tilde{Z}^{r.n} = BG \times \tilde{Z}^{r.n} \longrightarrow \tilde{Z}^{r.n}$$

qu'on peut interpréter comme une application équivariante

$$EG \times (Z^r)^{\otimes n} \longrightarrow \tilde{Z}^{r.n}$$

ou encore

$$EG \longrightarrow \mathrm{Hom}((Z^r)^{\otimes n}, \tilde{Z}^{r.n})$$

Dans le cas de la théorie $A^* = \mathcal{D}^*$, cette application est l'unique morphisme équivariant (à homotopie près) rendant commutatif le diagramme suivant :

$$\begin{array}{ccc} EG & \longrightarrow & \mathrm{Hom}((Z^r)^{\otimes n}, \tilde{Z}^{r.n}) \\ \downarrow & & \downarrow \\ k & \longrightarrow & \mathrm{Hom}((Z^r)^{\overline{\otimes} n}, \tilde{Z}^{r.n}) \end{array}$$

Introduisons maintenant un espace paramètre X et remplaçons les espaces du type $Z^m$ par les espaces $Z^m(X)$ qu'on peut interpréter comme les k-modules simpliciaux

$$s \mapsto Z^m(X \times \Delta_s)$$

Nous avons alors un diagramme commutatif

$$\begin{array}{ccc} EG & \longrightarrow & \mathrm{Hom}((Z^r)^{\otimes n}(X), \tilde{Z}^{r.n}(X)) \\ \downarrow & & \downarrow \\ k & \longrightarrow & \mathrm{Hom}((Z^r)^{\overline{\otimes} n}(X), \tilde{Z}^{r.n}(X)) \end{array}$$

**6.12.** Pour aller plus loin, nous allons d'abord montrer que le k-module simplicial $s \mapsto Z^m(X \times \Delta_s)$ (vu comme un complexe de chaînes par la correspondance de Dold-Kan) est quasi-isomorphe au complexe (homologique) suivant, noté $\mathcal{C}$ :

$$A^0(X) \longrightarrow \ \ldots \ \longrightarrow A^{m-1}(X) \longrightarrow Z^m(X) \longrightarrow 0$$

En effet, considérons le bicomplexe

$$\begin{array}{ccccccccc} 0 & & & & 0 & & 0 & & \\ \uparrow & & & & \uparrow & & \uparrow & & \\ A^0(X) & \longrightarrow & \ldots & \longrightarrow & A^{m-1}(X) & \longrightarrow & Z^m(X) & \longrightarrow & 0 \\ \uparrow & & & & \uparrow & & \uparrow & & \\ A^0(X \times \Delta_1) & \longrightarrow & \ldots & \longrightarrow & A^{m-1}(X \times \Delta_1) & \longrightarrow & Z^m(X \times \Delta_1) & \longrightarrow & 0 \\ \uparrow & & & & \uparrow & & \uparrow & & \\ A^0(X \times \Delta_2) & \longrightarrow & \ldots & \longrightarrow & A^{m-1}(X \times \Delta_2) & \longrightarrow & Z^m(X \times \Delta_2) & \longrightarrow & 0 \\ \uparrow & & & & \uparrow & & \uparrow & & \end{array}$$

Toutes les colonnes sont acycliques à l'exception éventuellement de la dernière (car les groupes d'homotopie des $A^j$ sont nuls). L'homologie $H_s$ du complexe total est ainsi isomorphe à $H_s(Z^m(X \times \Delta_*))$. Par ailleurs l'homologie des lignes ne change pas si on remplace les $X \times \Delta_j$ par X (lemme de Poincaré), auquel cas les lignes $C_{2i}$ et $C_{2i-1}$ deviennent isomorphes pour $i > 0$ par les flèches verticales. L'homologie du complexe total est donc aussi isomorphe à l'homologie de la première ligne. Ainsi, le complexe total associé à ce bicomplexe est quasi-isomorphe à la fois à sa première ligne et à sa dernière colonne : celles-ci sont donc homotopiquement équivalentes. Ce raisonnement montre par exemple que le conoyau de la flèche $Z^m(X \times \Delta_1) \longrightarrow Z^m(X)$ s'identifie à $H^m(X)$.
La même argumentation s'applique au produit tensoriel du k-module simplicial $s \mapsto Z^r(X \times \Delta_s)$, n fois par lui-même, compte tenu de la remarque suivante : il suffit de remplacer les lignes du bicomplexe précédent par des k-module simpliciaux en utilisant de nouveau le théorème de Dold-Kan. Ces lignes apparaissent ainsi comme les complexes normalisés de ces modules. En résumé, on obtient bien ainsi un quasi-isomorphisme (en zigzag) entre les complexes $(Z^r(X \times \Delta_s)^{\otimes n}$ et le produit tensoriel $\mathfrak{S}^{\otimes n}$, où $\mathfrak{S}$ est le complexe simplicial défini plus haut :

$$A^0(X) \longrightarrow \ \dots \ \longrightarrow A^{r-1}(X) \longrightarrow Z^r(X) \longrightarrow 0$$

Par ailleurs, $(Z^r(X \times \Delta_s)^{\otimes n}$ est aussi quasi-isomorphe au complexe $\mathfrak{S}^{\overline{\otimes} n}$. De cette manière, l'application E précédente :

$$EG \longrightarrow [\ s \mapsto \mathrm{Hom}(Z^{r.n}(X \times \Delta_s), \tilde{Z}^{r.n}(X \times \Delta_s))\ ]$$

s'interprète comme un morphisme équivariant de complexes bien défini à homotopie près (r étant choisi arbitrairement grand)

(E') $$C(EG) \longrightarrow \mathrm{Hom}(A^*(X)^{\otimes n}, A^*(X))$$

où C(EG) est la résolution acyclique libre canonique de k en tant que k[G]-module (notée précédemment EG par abus d'écriture). Cette application est donc bien un cup-produit supérieur (par adjonction) qui induit bien sur les cycles les opérations de Steenrod comme il a été indiqué plus haut.

## 7. Itération de la "bar-construction". Type d'homotopie.

**7.1.** Si X est un ensemble simplicial, nous désignons par $C_*(X\ ;\ k) = C_*(X)$ le groupe des chaînes sur X, à coefficients dans un groupe abélien arbitraire k.

Pour une fibration où B est muni d'un point base

$$F \longrightarrow E \longrightarrow B$$

nous avons une "augmentation" évidente

$$C_*(F) \longrightarrow \overline{\mathrm{Tot}}\ [C_*(E) \longrightarrow C_*(B \times E) \longrightarrow C_*(B \times B \times E) \longrightarrow \dots ]$$

Ici le complexe

$$C_*(E) \longrightarrow C_*(B \times E) \longrightarrow C_*(B \times B \times E) \longrightarrow \dots$$

est le double complexe de chaînes associé à l'ensemble cosimplicial

$$E \longrightarrow B \times E \longrightarrow B \times B \times E \longrightarrow \dots$$

introduit par D. Rector [20], dont l'importance a été remarquée dans un article de W.G. Dwyer ([2] p. 256). Le complexe total, noté $\overline{\mathrm{Tot}}$ , est défini en considérant le <u>produit</u> des éléments situés sur les diagonales appropriées (cf. [9] p. 8). Nous donnerons plus de précision sur le degré total un peu plus loin.

**7.2.** D'après [5] (voir aussi l'annexe écrite par Michel Zisman et l'article de Brooke Shipley cité dans l'introduction en Note 3), cette application d'augmentation est un quasi-isomorphisme si le groupe fondamental $\pi_1(B, b)$ de B opère de manière nilpotente sur $H_*(F ; k)$, b étant le point base de B (cependant, B n'est pas nécessairement connexe, contrairement à ce qui est supposé dans [5]). Ceci a lieu automatiquement dans deux cas importants :

1. $\pi_1(B, b)$ est trivial, ou
2. $\pi_1(B, b)$ est un p-groupe fini et $k = \mathbf{Z}/p^\alpha$ pour un certain $\alpha$.

Si E est contractile, F est homotopiquement équivalent à l'espace des lacets pointé de B noté $\Omega(B)$. Par conséquent, nous avons un quasi-isomorphisme

$$C_*(\Omega(B)) \approx \overline{\mathrm{Tot}}\,[\, k \longrightarrow C_*(B) \longrightarrow C_*(B \times B) \longrightarrow \dots ]$$

**7.3.** De manière plus précise, si B est un ensemble pointé, de point base $*$, on peut lui associer un **ensemble cosimplicial** $L^*B$, défini par $L^nB = B^n$ , avec $B^0 = *$, les premiers opérateurs de coface

$$* \xrightarrow{d_0,\, d_1} B^1 \xrightarrow{d_0, d_1, d_2} B^2 \xrightarrow{d_0, d_1, d_2, d_3} B^3$$

étant définis ainsi

$d_0(*) = d_1(*) = *$

$d_0(x) = (*, x), d_1(x) = (x, x), d_2(x) = (x, *)$

$d_0(x, y) = (*, x, y), d_1(x, y) = (x, x, y), d_2(x, y) = (x, y, y), d_3(x, y) = (x, y, *)$

**7.4.** Au niveau des chaînes, ceci revient à considérer le bicomplexe concentré sur le deuxième quadrant et défini par $C_{r,s}(B) = C_{-r}(B^s)$, $C_*$ désignant le complexe des chaînes singulières, les différentielles étant induites par les $d_i$ et la différentielle usuelle sur les chaînes

$$\begin{array}{ccccc}
\downarrow & & \downarrow & & \downarrow \\
C_2(B^2) & \longleftarrow & C_2(B) & \longleftarrow & C_2(*) = k \\
\downarrow & & \downarrow & & \downarrow \\
C_1(B^2) & \longleftarrow & C_1(B) & \longleftarrow & C_1(*) = k \\
\downarrow & & \downarrow & & \downarrow \\
C_0(B^2) & \longleftarrow & C_0(B) & \longleftarrow & C_0(*) = k
\end{array}$$

Le “complexe total” $\overline{\mathrm{Tot}}$ défini par

$$\overline{\mathrm{Tot}}_n = \prod_{u - v = n} C_u(B^v)$$

dont l’homologie calcule celle de l’espace des lacets $\Omega(B)$ (ce résultat classique est dû à Adams et Hilton).

**7.5.** Nous voulons maintenant itérer cette procédure. Ceci peut se faire par exemple si B est 2-connexe (donc $\Omega(B)$ simplement connexe) ou si les deux premiers groupes d’homotopie sont des p-groupes finis (avec $k = \mathbf{Z}/p^{\alpha}$). Si nous remplaçons $C_*(B^n)$ par

$$\overline{\mathrm{Tot}}\,[k \longrightarrow C_*(B^n) \longrightarrow C_*(B^n \times B^n) \longrightarrow ...]$$

nous trouvons l’isomorphisme suivant dans la catégorie dérivée :

$$C_*(\Omega^2(B)) \approx \overline{\mathrm{Tot}}[\, C_*(B^{mn})]$$

où par convention $B^0$ est un point et où $\overline{\mathrm{Tot}}$ signifie encore le complexe total (complété) associé à un triple complexe. Par récurrence sur r, on démontre ainsi le théorème suivant :

**7.6. THEOREME.** *Supposons que* B *soit* r-*connexe ou que les* r *premiers groupes d’homotopie de* B *soient des* p-*groupes finis (auquel cas nous supposons en outre que* $k = \mathbf{Z}/p^{\alpha}$). *Nous avons alors un isomorphisme dans la catégorie dérivée* :

$$C_*(\Omega^r(B)) \approx \overline{\mathrm{Tot}}[\, C_*(B^{v_1 \cdots v_r})]$$

$\overline{\mathrm{Tot}}$ *désignant le complexe total associé au* (r+1) -*complexe évident.*

**7.7. Remarque.** En considérant la diagonale dans le k-module cosimplicial $C_*(B^{v_1 \cdots v_r})$, on peut donner une définition homotopiquement équivalente de $\overline{\mathrm{Tot}}$, analogue à celle décrite en 7.4 et qui est la suivante

$$\overline{\mathrm{Tot}}_n = \prod_{u - v = n} C_u(B^{v^r})$$

**7.8.** Nous pouvons traduire ces considérations dans un contexte cohomologique. Nous avons alors une application d’augmentation pour les <u>cochaînes</u> :

$$\mathrm{Tot}\,[.... \longrightarrow C^*(B \times B \times E) \longrightarrow C^*(B \times E) \longrightarrow C^*(E)] \longrightarrow C^*(F)$$

où “Tot” désigne le complexe total classique : il est obtenu en considérant les **sommes** (au lieu des produits) d’éléments situés sur les diagonales appropriées. Si nous supposons que k est un corps et que les homologies de E, B et F sont des espaces vectoriels de dimension finie en chaque dimension, il est facile de montrer que la cohomologie de ce complexe total est duale de l’homologie du complexe $\overline{\mathrm{Tot}}$ précédent. Il en résulte que l’application d’augmentation est aussi un quasi-isomorphisme en cohomologie. Ceci se produit par exemple quand $k = \mathbf{Z}/p\mathbf{Z}$ et quand tous les groupes d’homotopie des espaces impliqués sont des p-groupes finis. D'après le théorème précédent, nous avons par conséquent un isomorphisme naturel dans la catégorie dérivée entre $C^*(\Omega^r(B))$ and $Tot[C^*(B^{n_1 \dots n_r})]$. En appliquant le lemme des cinq, on en déduit finalement le théorème suivant :

**7.9. THEOREME.** *Supposons que tous les groupes d'homotopie de* B *soient des* p-*groupes finis. Nous avons alors l'isomorphisme naturel suivant dans la catégorie dérivée* :

$$C^*(\Omega^r(B) ; \mathbf{Z}/p^\alpha\mathbf{Z}) \approx Tot[C^*(B^{n_1 \dots n_r}) ; \mathbf{Z}/p\mathbf{Z}]$$

Plus précisément, en considérant de nouveau la diagonale, le complexe $C^*(\Omega^r(B) ; \mathbf{Z}/p\mathbf{Z})$ est quasi-isomorphe au complexe

$$\overline{Tot}^n = \bigoplus_{u - v = n} C^u(B^{v^r})$$

**7.10. Remarques importantes.** En changeant le signe de v, il sera commode par la suite de considérer le double complexe cohomologique situé sur le quatrième quadrant défini par

$$C^{u,w} = \bigoplus_{u + w = n} C^u(B^{(-w)^r})$$

On notera que les différentielles vont bien dans le bon sens, soit de $C^{u,w}$ vers $C^{u+1,w}$ ou $C^{u,w+1}$. Par ailleurs, si B est simplement connexe, nous pouvons généraliser ce théorème : il suffit de supposer que les r premiers groupes d'homotopie sont des p-groupes finis et que $H_i(B ; \mathbf{Z}/p)$ est de dimension finie pour chaque i. Si ces r groupes d'homotopie sont nuls, le corps k peut être quelconque, mais on doit toujours supposer que $H_i(B ; k)$ est de dimension finie pour tout i.

Enfin, dans les considérations précédentes, nous pouvons remplacer le complexe classique des cochaînes $C^*(B^m ; k)$ par le complexe des formes différentielles quasi-commutatives $\mathcal{D}^*(B^m)$, défini dans les paragraphes précédents puisque ces complexes sont naturellement quasi-isomorphes.

**7.11.** L'homologie de Hochschild r-itérée ${}_{(r)}H_*(A)$ d'une algèbre **commutative**[19] augmentée A est définie ainsi : à l'algèbre A nous associons le groupe abélien multisimplicial

$$(m_1, \dots, m_r) \mapsto A^{\otimes m_1 \dots m_r}$$

Les opérations face et dégénérescence sont évidentes sur chaque "coordonnée". On a par exemple (ε désignant l'augmentation)

$$\partial_0(a_1, \dots, a_m) = (\varepsilon(a_1)a_2, \dots, a_m),$$
$$\partial_1(a_1, \dots, a_m) = (a_1a_2, \dots, a_m),$$
$$\partial_m(a_1, \dots, a_m) = (a_1, \dots, a_{m-1}\varepsilon(a_m))$$

L'homologie totale $H_*(\mathcal{C})$ du complexe $\mathcal{C}$ obtenu est cette r-itération de l'homologie de Hochschild. Si r = 1, nous retrouvons l'homologie de Hochschild usuelle à coefficients dans k (via l'augmentation). Si r > 1, cette homologie a été étudiée de manière intensive par T. Pirashvili[20] dans le cas où A est une algèbre commutative sur un corps k de caractéristique

[19] ou commutative au sens gradué si A est graduée.

0. Cette homologie itérée fait partie d'un arsenal bien connu d'algèbre homologique qui est l'itération de la bar construction d'une algèbre augmentée simpliciale <u>commutative</u>. De manière précise, considérons une algèbre de ce type que nous écrivons sous la forme

$$\longrightarrow A_n \longrightarrow A_{n-1} \longrightarrow \ldots \longrightarrow A_1 \longrightarrow k$$

où k est concentré en degré (simplicial) 0 et $A_n$ en degré n. On suppose bien entendu que les opérateurs face et dégénérescence sont des morphismes d'algèbres commutatives. La bar-construction appliquée à cette algèbre donne naissance à une algèbre <u>bisimpliciale</u> commutative

$$\begin{array}{ccccccccc}
\longrightarrow & (A_3)^{\otimes 3} & \longrightarrow & (A_2)^{\otimes 3} & \longrightarrow & (A_1)^{\otimes 3} & \longrightarrow & k \\
 & \downarrow & & \downarrow & & \downarrow & & \downarrow \\
\longrightarrow & (A_3)^{\otimes 2} & \longrightarrow & (A_2)^{\otimes 2} & \longrightarrow & (A_1)^{\otimes 2} & \longrightarrow & k \\
 & \downarrow & & \downarrow & & \downarrow & & \downarrow \\
\longrightarrow & (A_3)^{\otimes 1} & \longrightarrow & (A_2)^{\otimes 1} & \longrightarrow & (A_1)^{\otimes 1} & \longrightarrow & k \\
 & \downarrow & & \downarrow & & \downarrow & & \downarrow \\
\longrightarrow & k & \longrightarrow & k & \longrightarrow & k & \longrightarrow & k
\end{array}$$

dont la diagonale est simplement l'algèbre simpliciale commutative suivante

$$\longrightarrow (A_3)^{\otimes 3} \longrightarrow (A_2)^{\otimes 2} \longrightarrow A_1 \longrightarrow k$$

qui peut être itérée de nouveau. Ainsi, la 2e itération de l'algèbre commutative de départ A est l'algèbre simpliciale

$$\longrightarrow (A)^{\otimes 9} \longrightarrow (A)^{\otimes 4} \longrightarrow A \longrightarrow k$$

tandis que la 3e itération est

$$\longrightarrow (A)^{\otimes 27} \longrightarrow (A)^{\otimes 8} \longrightarrow A \longrightarrow k$$

**7.12.** Une théorie analogue peut être développée plus généralement pour les ADG <u>quasi-commutatives</u>. Dans les définitions précédentes, pour $A = \mathcal{D}^*(B)$ par exemple, il suffit de remplacer $A^{\otimes m_1 \ldots m_r}$ par $\mathcal{D}^*(B)^{\overline{\otimes} m_1 \ldots m_r}$. D'après nos motivations topologiques, les complexes totaux (situés sur le quatrième quadrant d'après 7.10) doivent être considérés dans le sens classique, en faisant la <u>somme</u> d'éléments situés sur les diagonales appropriées.

**7.13.** De manière plus précise, limitons nous provisoirement à r = 2 pour mieux illustrer notre raisonnement. Le morphisme de face

$$b_i : (A^{\otimes n})^{\otimes p} \longrightarrow (A^{\otimes n})^{\otimes (p-1)}$$

par exemple n'est pas un homomorphisme d'anneaux, ce qui empêche la commutativité des diagrammes "cubiques" qu'il faudrait écrire dans une bar-construction "tridimensionnelle". Bien que les diagrammes

[20] Annales Scientifiques de l'Ecole Normale Supérieure (33) 2000, p. 151-179.

$$\begin{array}{ccc} A^{\otimes np} & \longrightarrow & A^{\otimes (n-1)p} \\ \downarrow & & \downarrow \\ A^{\otimes n(p-1)} & \longrightarrow & A^{\otimes (n-1)(p-1)} \end{array}$$

obtenus en appliquant des opérateurs face $b_i$ successifs, ne soient pas commutatifs en général, ils le sont "à homotopie près" dans le sens suivant . **Il convient de remplacer tous les produits tensoriels par des produits tensoriels réduits**. Si nous désignons par $A^{\overline{\otimes} r}$ le k-module gradué $\mathcal{D}^*(B)^{\otimes r}$, nous obtenons bien alors des diagrammes commutatifs

$$\begin{array}{ccc} A^{\overline{\otimes} np} & \longrightarrow & A^{\overline{\otimes} (n-1)p} \\ \downarrow & & \downarrow \\ A^{\overline{\otimes} n(p-1)} & \longrightarrow & A^{\overline{\otimes} (n-1)(p-1)} \end{array}$$

**7.14.** Le bicomplexe $C^{u,w}$ considéré en 7.10 peut s'écrire[21] ainsi de manière équivalente[22] en remplaçant $C^{u,w} = \bigoplus_{u+w=n} C^u(B^{(-w)^r}$ par

$$D^{u,w} = \bigoplus_{u+w=n} [\, \mathcal{D}^*(B^{\overline{\otimes}(-w)^r}]^u$$

où le symbole [ ]$^u$ signifie qu'on se restreint aux éléments de degré u. Considérons maintenant la suite spectrale associée à ce double complexe. Le terme $E_1$ n'est autre que

$$\bigoplus_{u+w=n} [\, H^*(B^{\overline{\otimes}(-w)^r}]^u$$

Le terme $E_2$ est simplement l'homologie de Hochschild r-itérée de l'algèbre commutative (au sens gradué) H*(B) = C. De manière précise, ${}_{(r)}H_t(C)$ hérité d'une graduation (cohomologique) induite par celle de C. Si on note ${}_{(r)}H_t(C)^u$ le terme de degré u pour cette graduation, on voit finalement que le terme $E_2$ de la suite spectrale est défini par

$$E_2^{u,w} = {}_{(r)}H_t(C)^u$$

avec t = -w. D'après les considérations topologiques précédentes, cette suite spectrale converge vers la cohomologie de degré u + w de l'espace de lacets r-itéré de B (voir aussi l'article de Brooke Shipley cité dans l'introduction en note 3).

**7.14.** Un cas particulièrement intéressant est le calcul de $H^0(\Omega^r B) = \mathrm{Hom}_{Ens}(\pi_r(B), k)$ qui serait l'aboutissement des termes $E_2$ correspondant aux groupes d'homologie de Hochschild itérés ${}_{(r)}H_u(C)^u$. Une majoration de la dimension de la somme de ces espaces vectoriels de cohomologie donnerait naissance à une majoration de l'ordre du groupe d'homotopie $\pi_r(B)$. Malheureusement, nos connaissances actuelles sur ces groupes d'homologie de Hochschild

[21] A partir de maintenant, on suppose que k est le corps **Z**/p**Z**, ce qui permet d'appliquer la formule de Künneth.

[22] avec $u \geq 0$ et $w \leq 0$.

itérés restent encore très fragmentaires.

**7.15.** Nous n'avons pas encore tenu compte d'une donnée fondamentale de la bar-construction qui est l'existence d'un coproduit. Celui-ci repose sur le lemme suivant (que nous écrivons simplement pour la double itération, mais qui peut se transcrire aisément dans le cas général).

**7.16. LEMME.** *L'homomorphisme évident*

$$\bigoplus_{(m,n)} A^{\otimes mn} \longrightarrow \bigoplus_{(m,\,n,\,r)} A^{\otimes (m-r)n} \otimes A^{\otimes rn}$$

*induit un morphisme sur les produits tensoriels réduits correspondants*

$$\phi : \bigoplus_{(m,n)} A^{\overline{\otimes} mn} \longrightarrow \bigoplus_{(m,\,n,\,r)} A^{\overline{\otimes}(m-r)n} \otimes A^{\overline{\otimes} rn}$$

*Il munit ainsi la somme directe (pour tout n) des* $A^{\overline{\otimes} mn}$ *d'une structure de cogèbre compatible avec les différentielles du bicomplexe.*

*Démonstration.* Posons

$S = \{0, ..., m\} \times \{0, ..., n\}$

$S_r = \{0, ..., m-r\} \times \{0, ..., n\}$

$S'_r = \{m-r+1, ..., m\} \times \{0, ..., n\}$.

La décomposition $S = S_r \cup S'_r$ induit des morphismes $A^{\otimes S} \longrightarrow A^{\otimes S_r} \otimes A^{\otimes S'_r}$ et $A^{\overline{\otimes} S} \longrightarrow A^{\overline{\otimes} S_r} \otimes A^{\overline{\otimes} S'_r}$ avec des notations évidentes. Les opérateurs face du bicomplexe sont associés à des morphismes "rectangulaires" $S \longrightarrow T$ qui sont compatibles avec le morphisme $\phi$ en raison de la remarque générale faite en 1.1 : ces opérateurs face ne dépendent pas du choix d'un ordre sur S ou sur T.

En généralisant ces considérations à la r-itération, on voit par exemple que le terme $E_2$ de la suite spectrale décrite en 7.13 et 7.14 hérite d'une structure de cogèbre qui devrait être compatible (sur l'aboutissement) avec la structure de cogèbre de $H^*(\Omega^r B)$, induite par la structure d'espace de lacets de $\Omega^r B$.

Enonçons maintenant le résultat fondamental de ce paragraphe, dont les résultats précédents sont annonciateurs :

**7.17. THEOREME.** *Considérons deux ensembles simpliciaux* X *et* Y *connexes, nilpotents et* p*-complets de type fini* [23]. *Nous supposons qu'il existe une suite en zigzag de quasi-isomorphismes d'ADG quasi-commutatives* (avec $k = \mathbf{F}_p$)

$$\mathcal{D}^*(X) \longrightarrow A \longleftarrow B \longrightarrow .... \longleftarrow \mathcal{D}^*(Y)$$

*Alors* X *et* Y *ont le même type d'homotopie.*

*Démonstration.* Soit Z un ensemble simplicial. D'après [10], il existe un quasi-isomorphisme naturel entre les k-modules différentiels gradués $\mathcal{D}^*(Z)$ et $C^*(Z)$. Par ailleurs, nous pouvons associer à une ADGQ une $E_\infty$-algebre en utilisant la méthode de Kriz et May [13]. Plus précisément, une ADGQ est un cas particulier de ce que ces auteurs appellent une “algèbre partielle” (cf. [13], p. 40). Soit $\mathcal{P}$ (resp. $\mathcal{E}$, resp $\mathcal{ESP}$) la catégorie des algèbres partielles (resp. des $E_\infty$-algèbres, resp. des $E_\infty$-algèbres partielles simpliciales). Dans [13] on décrit un diagramme de catégories et de foncteurs qui est commutatif à isomorphisme près ($\varphi$ et $\psi$ étant des quasi-isomorphismes des modules différentiels gradués sous-jacents)

$$\begin{array}{ccc} & & \mathcal{P} \\ & \overset{Id}{\nearrow} & \uparrow \varphi \\ \mathcal{P} & \xrightarrow{\ V\ } & \mathcal{ESP} \\ & \underset{W}{\searrow} & \downarrow \psi \\ & & \mathcal{E} \end{array}$$

[23] Ceci veut dire que sa tour de Postnikov peut être choisie en sorte que chaque fibre est de type $K(\mathbf{Z}/p, n)$ ou $K(\hat{\mathbf{Z}}_p, n)$, $\hat{\mathbf{Z}}_p$ désignant l'anneau des entiers p-adiques [8].

Les quasi-isomorphismes figurant dans nos hypothèses

$$\mathcal{D}^*(X) \longrightarrow A \longleftarrow B \longrightarrow .... \longleftarrow \mathcal{D}^*(Y)$$

impliquent par conséquent une suite de quasi-isomorphismes entre les $E_\infty$-algèbres associées via le foncteur W.

Par ailleurs, d’après un résultat récent de Mandell [20], pour chaque ensemble simplicial Z, les $E_\infty$-algèbres $\mathcal{D}^*(Z)$ et $C^*(Z)$ sont aussi reliées par une suite de quasi-isomorphismes de $E_\infty$-algèbres. D’après la conclusion précédente, nous en déduisons que $C^*(X)$ et $C^*(Y)$ sont aussi reliés par une suite de quasi-isomorphismes de $E_\infty$-algèbres.

Puisque X et Y sont nilpotents, un deuxième résultat clé de Mandell [19] implique que X et Y ont le même type d’homotopie, ce qui conclut la démonstration de notre théorème. Une version plus faible du théorème est la suivante :

**7.18. THEOREME.** *Considérons deux ensembles simpliciaux* <u>*finis*</u> *et connexes* X *et* Y *tels que leurs groupes d’homotopie soient des* p-*groupes finis. Nous supposons qu’il existe une suite en zigzag de quasi-isomorphismes d’algèbres différentielles graduées quasi-commutatives* (avec $k = \mathbf{F}_p$)

$$\mathcal{D}^*(X) \longrightarrow A \longleftarrow B \longrightarrow .... \longleftarrow \mathcal{D}^*(Y)$$

*Alors* X *et* Y *ont le même type d’homotopie.*

**7.19.** En fait, nous avons décrit dans les paragraphes précédents une procédure explicite pour déterminer algébriquement les groupes d’homotopie de X, via une “itération” convenable de la bar-construction [5], à partir de l’ADG quasi-commutative $A = \mathcal{D}^*(X)$. Plus précisément, la correspondance $(m_1,\dots, m_r) \mapsto A_{m_1\dots m_r}$ définit un module gradué r-simplicial (le point base est utilisé pour définir certains opérateurs face) dont la cohomologie est celle de l’espace de lacets r-itéré de X, noté ici $\Omega^r(X)$.

**7.20. Remarques**. Récemment, Mandell a annoncé une extension de son résultat principal à k = **Z**. Il en résulte aussitôt une extension du théorème 7.18 à ce cas. Il serait évidemment intéressant de démontrer le théorème 7.18 de manière indépendante sans faire le détour par les $E_\infty$-algèbres. Il reste aussi le problème ouvert d'un modèle "minimal" dans la théorie des ADGQ, tel qu'il existe dans celles des ADG commutatives par la théorie de Sullivan [28].

## 8. Autres structures algébrique sur $\mathcal{D}^*(X)$ : tressage, opérateur de translation, produit réduit...

Nous commencerons par quelques généralités sur les algèbres tressées que nous détaillons ici pour la commodité du lecteur.

**8.1.** Soient k un anneau commutatif quelconque et A un k-module. Un "tressage" sur A est la donnée d'un k-**automorphisme**

$$R : A \otimes A \longrightarrow A \otimes A$$

(noté aussi parfois $R_A$) satisfaisant aux "équations de Yang-Baxter" :

$$R_{12}.R_{23}.R_{12} = R_{23}.R_{12}.R_{23}$$

avec les notations standard pour les endomorphismes de $A \otimes A \otimes A$ (cf. [15] par exemple).

**8.2. DEFINITION.** *Soit* A *une algèbre munie d'un tressage* R *en tant que* k-*module. Alors* A *est dite tressée en tant qu'algèbre si*

1. $R(1 \otimes a) = a \otimes 1$ *et* $R(a \otimes 1) = 1 \otimes a$
2. *Nous avons les identités* $\mu_{23}.R_{12}.R_{23} = R.\mu_{12}$ *et* $\mu_{12}.R_{23}.R_{12} = R.\mu_{23}$
3. *Le diagramme suivant commute*

$$\begin{array}{ccc} A \otimes A & \xrightarrow{R} & A \otimes A \\ & \searrow \quad \swarrow & \\ & A & \end{array}$$

*les flèches obliques représentant la multiplication dans l'algèbre.*
*En outre, si A est une algèbre différentielle graduée, nous supposerons que le tressage R est de degré 0 et commute aux différentielles évidentes.*

**8.3. Exemples**.

**8.3.1.** Toute algèbre commutative est évidemment tressée commutative, l'endomorphisme R étant défini par la formule $R(a \otimes b) = b \otimes a$. De manière plus générale, soit A une algèbre graduée commutative (au sens gradué). Alors A est aussi tressée commutative, le tressage étant défini sur des éléments homogènes a et b par la formule suivante :

$$R(a \otimes b) = (-1)^{\deg(a)\deg(b)}\, b \otimes a$$

**8.3.2.** Comme l'a remarqué P. Nuss[24], <u>toute</u> algèbre A peut être tressée en posant

$$a \otimes b \mapsto ab \otimes 1 - a \otimes b + 1 \otimes ab$$

Ce tressage donne lieu à une représentation du groupe des tresses $\mathfrak{B}_n$ dans le groupe des automorphismes de $A^{\otimes n}$, représentation qui se factorise en fait par le groupe symétrique $\mathfrak{S}_n$.

**8.3.3.** Soit A une algèbre de Hopf. Alors S. L. Woronowicz[25] a défini un tressage intéressant sur A (non involutif en général). Par exemple, si A = k[G] est l'algèbre d'un groupe discret G, ce tressage est simplement induit par

$$g \otimes h \mapsto h \otimes h^{-1}gh$$

pour des éléments g et h de G. Voici un autre exemple : si $A = U(\mathfrak{g})$ est l'algèbre enveloppante d'une algèbre de Lie $\mathfrak{g}$, le tressage est induit par

$$g \otimes h \mapsto h \otimes g + 1 \otimes [g, h]$$

**8.3.4.** Voici maintenant l'exemple qui nous intéresse le plus dans cet article (les vérifications techniques des axiomes sont faciles, mais fastidieuses). Dans [14] nous donnons une interprétation plus conceptuelle de ces calculs dans un contexte plus général.

Considérons une k-algèbre **commutative** $\Lambda$ munie d'un automorphisme d'algèbre $\alpha$ noté aussi $a \mapsto \bar{a} = \alpha(a)$ ($\alpha$ n'est pas nécessairement involutif). Soit $A = \Lambda \oplus \Omega^1(\Lambda)$, où $\Omega^1(\Lambda)$ désigne le k-module des différentielles de Käehler "tordu" par cet automorphisme : c'est le quotient du k-module des formes différentielles non commutatives de degré 1 par le sous-module engendré par les relations suivantes :

$$u.d(vv') = uv.dv' + u.dv.v'$$

$$u.dv.v' = u\overline{v'}.dv$$

Le k-module $\Omega^1(\Lambda)$ peut être aussi défini comme le conoyau de $\bar{b}$, où

$$\bar{b} : \Lambda^{\otimes 3} \longrightarrow \Lambda^{\otimes 2}$$

est l'opérateur bord de Hochschild "tordu", soit

$$\bar{b}(a_0 \otimes a_1 \otimes a_2) = a_0a_1 \otimes a_2 - a_0 \otimes a_1a_2 + \bar{a}_2a_0 \otimes a_1$$

Posons $\Omega^0(\Lambda) = \Lambda$ et $\Omega^i(\Lambda) = 0$ pour $i > 1$. Nous voyons alors que $A = \bigoplus_i \Omega^i(\Lambda)$ est une algèbre différentielle graduée de manière évidente. Un tressage R sur A est défini par les formules suivantes

1) $R(u \otimes v) = v \otimes u$ si u et v sont de degré 0.

2) $R(udv \otimes w) = \overline{w} \otimes udv$ , pour udv de degré 1 et w de degré 0. De même,

[24] Noncommutative descent and non-abelian cohomology. K-Theory 12 (1997), N°. 1, 23--74.

[25] Solutions of the braid equation related to a Hopf algebra. Lett. Math. Phys. 23 (1991)

3) $R(w \otimes udv) = udv \otimes w + u(v - \bar{v}) \otimes dw$

4) $R(udv \otimes wdt) = - \overline{wdt} \otimes udv$ , pour udv et wdt de degré 1.

Bien entendu, il convient de vérifier que ces formules ont un sens en utilisant les identités du type $u.d(vv') = uv.dv' + u\overline{v'}.dv$, ce qui ne présente pas de difficultés particulières.

**Remarque** : plus généralement, on peut définir de même des formes différentielles “tordues” de degrés > 1 : cf. l'article [14] déjà cité plus haut.

**8.3.4.1.** Deux exemples illustrant la théorie esquissée en 8.3.4 méritent d'être mentionnés. Le premier est bien connu [18] et consiste à choisir pour $\Lambda$ l'algèbre k[t] des polynômes à une variable t et pour $\varphi$ l'automorphisme $t \mapsto qt$, où q est un élément inversible de k. On obtient alors l'algèbre des formes différentielles polynomiales “quantiques” à une variable. Par exemple, la différentielle du monôme $t^n$ est égale à $[n]_q\, t^{n-1}dt$ , où $[n]_q$ est l'entier “quantique” $\dfrac{q^n - 1}{q-1}$. Un autre exemple plus élaboré et utilisant des séries convergentes est décrit dans [12]. On peut remarquer que si $[n]_q$ est inversible dans k, le lemme de Poincaré est vrai : la cohomologie de l'algèbre différentielle graduée $A = \Lambda \oplus \Omega^1(\Lambda)$ est concentrée en degré 0 et est isomorphe à k.

**8.3.4.2.** Le deuxième exemple a été explicitement décrit dans [13], [14] et abondamment traité dans les paragraphes précédents. Nous le reprenons ici pour la commodité du lecteur, en insistant sur la notion de tressage.

Soit $\Lambda$ la k-algèbre $\mathcal{D}^0(x)$ des fonctions $f = f(x) : \mathbf{Z} \longrightarrow k$ qui sont constantes quand x tend vers $+ \infty$ ou $- \infty$ (deux limites indépendantes). Cette algèbre peut être munie de l'automorphisme $f \mapsto \bar{f}$, où $\bar{f}(x) = f(x+1)$ Le module des différentielles “tordues” de $\Lambda$ se décrit maintenant comme le quotient du module des différentielles non commutatives $\Omega^1_{nc}(\Lambda)$ par le sous-module engendré par les relations $(f.dg).h - (\bar{h}.f).dg$. Comme il est montré dans [13] et [14], ce module (noté $\mathcal{D}^1(x)$) s'identifie au sous $\Lambda$-module à gauche de $\mathcal{D}^0(x)$ formé des fonctions $\omega$ telles que $\omega(\infty) = \omega(-\infty) = 0$. Nous munissons $\mathcal{D}^1(x)$ d'une structure de $\Lambda$-module à droite par la formule suivante :

$$(\omega.f)(x) = (\bar{f}.\omega)(x) = f(x + 1).\omega(x)$$

avec $f \in \mathcal{D}^0(x)$. On montre alors que le module des différentielles tordues s'identifie à $\mathcal{D}^1(x)$ par l'application $f.dg \mapsto f.\bar{g} - f.g$ et que la différentielle

$$d : \mathcal{D}^0(x) \longrightarrow \mathcal{D}^1(x)$$

associe à une fonction f la “forme différentielle” $\omega(x) = f(x + 1) - f(x)$, ce qui n'est autre que le calcul aux différences.

D'après la théorie générale décrite plus haut, le tressage sur $\mathcal{D}^*(x)$ s'explicite ainsi (f et g étant de degré 0, $\omega$ et $\theta$ de degré 1) :

$R(f \otimes g) = g \otimes f$

$R(\omega \otimes g) = \bar{g} \otimes \omega$

$R(h \otimes f.dg) = f.dg \otimes h + f.(g - \bar{g}) \otimes dh$

$R(\omega \otimes \theta) = - \bar{\theta} \otimes \omega$

[l'automorphisme $\theta \mapsto \bar{\theta}$ étant induit par l'automorphisme $f \mapsto \bar{f}$ de $\Lambda$]

**8.4.** Ces généralités s'appliquent aussi au produit tensoriel de plusieurs facteurs. On en déduit que l'algèbre $\mathcal{D}^*(x_0,..., x_n) = \mathcal{D}^*(x_0) \otimes ... \otimes \mathcal{D}^*(x_n)$ est une ADG tressée. Plus généralement, en considérant r et s ≤ n, on a aussi un "tressage"

$$\mathcal{D}^*(x_0,..., x_r) \otimes \mathcal{D}^*(x_0,..., x_s) \longrightarrow \mathcal{D}^*(x_0,..., x_s) \otimes \mathcal{D}^*(x_0,..., x_r)$$

obtenu en plongeant les deux algèbres dans $\mathcal{D}^*(x_0,..., x_n)$.

**8.5.** Il nous reste à construire un" tressage" dans un sens généralisé évident et où nous supposons ici que X et Y sont des **complexes simpliciaux finis** comme dans le § 4 :

$$R_{X,Y} : \mathcal{D}^*(X) \otimes \mathcal{D}^*(Y) \longrightarrow \mathcal{D}^*(Y) \otimes \mathcal{D}^*(X)$$

Pour cela, on définit d'abord un "tressage bisimplicial"

$$\mathcal{D}^*(\Delta_\natural) \otimes \mathcal{D}^*(\Delta_\#) \xrightarrow{R} \mathcal{D}^*(\Delta_\#) \otimes \mathcal{D}^*(\Delta_\natural)$$

en écrivant simplement $\mathcal{D}^*(\Delta_r)$ (resp. $\mathcal{D}^*(\Delta_s)$) comme le quotient de $\mathcal{D}^*(x_0,..., x_r)$ (resp. $\mathcal{D}^*(x_0,..., x_s)$) et en utilisant le tressage

$$\mathcal{D}^*(x_0,..., x_r) \otimes \mathcal{D}^*(x_0,..., x_s) \longrightarrow \mathcal{D}^*(x_0,..., x_s) \otimes \mathcal{D}^*(x_0,..., x_r)$$

vu plus haut. En effet, comme on le voit par les formules, si $I_r$ désigne le noyau de l'application canonique de $\mathcal{D}^*(x_0,..., x_r)$ dans $\mathcal{D}^*(\Delta_r)$, l'image de $I_r \otimes \mathcal{D}^*(x_0,..., x_s)$ par le tressage est contenue dans

$$\mathcal{D}^*(x_0,..., x_s) \otimes I_s + I_r \otimes \mathcal{D}^*(x_0,..., x_s)$$

avec des notations évidentes. Une remarque analogue s'applique en permutant les deux facteurs. Par ailleurs, ce tressage commute aux opérations face (poser $x_i = \pm \infty$) [26].

**8.6.** Supposons maintenant que X et Y soient des **complexes** simpliciaux finis, sous-complexes de $\Delta_n$ par exemple (comme dans le § 4). Puisque le tressage sur les simplexes est compatible avec les opérations face, il suffit de remarquer que les tressages "locaux"

$$\mathcal{D}^*(P) \otimes \mathcal{D}^*(Q) \longrightarrow \mathcal{D}^*(Q) \otimes \mathcal{D}^*(P)$$

se globalisent sur X et Y (X étant une réunion de simplexes P et Y étant une réunion de simplexes Q). On en déduit le tressage annoncé

[26] Attention : ce tressage est compatible avec les opérations face, mais pas avec les opérations de dégénérescence.

$$\mathcal{D}^*(X) \otimes \mathcal{D}^*(Y) \longrightarrow \mathcal{D}^*(Y) \otimes \mathcal{D}^*(X)$$

Des considérations précédentes (en faisant X = Y), il résulte ainsi que $\mathcal{D}^*(X)$ est une ADG qui est non seulement quasi-commutative, mais aussi tressée. En particulier, sur chaque complexe simplicial fini, on a une représentation du groupe des tresses $\mathfrak{B}_n$ sur $\mathcal{D}^*(X)^{\otimes n}$. Plus généralement, si $\mathcal{F}$ est un faisceau en k-algèbres commutatives de base X et si $\mathcal{U}$ est un recouvrement ouvert **fini** de X, l'algèbre $\mathcal{D}^*(\mathcal{U}\ ;\ \mathcal{F}) = C^*(\mathcal{U}^{\natural}\ ;\ \mathcal{F})\ \nabla\ \mathcal{D}^*(\Delta_{\natural})$ est une algèbre tressée.

Voici maintenant une autre structure qui présente aussi un certain intérêt et qui est induite par la translation de **Z**. En effet, celle-ci induit un automorphisme T de l'algèbre $\mathcal{D}^*(x)$ qui respecte les augmentations obtenues en posant $x = \pm\ \infty$.

**8.7. Lemme.** *L'automorphisme* T *sur* l'ADG $\mathcal{D}^*(x)$ *est homotope à l'identité.*
*Démonstration.* On définit un opérateur d'homotopie K de degré -1. Il est 0 sur les fonctions de Heaviside et la fonction unité. Il associe à la fonction de Dirac $\omega_x$ (vue comme forme différentielle de degré un) associe la même fonction de Dirac (notée $\delta_x$ ), vue maintenant comme élément de degré 0. Il est clair qu'on a alors dK + Kd = T - 1.

**8.8.** Il en résulte que l'opérateur $T \otimes T \otimes ... \otimes T$ (n+1 facteurs) est aussi homotope à $1 \otimes 1 \otimes ... \otimes 1$ . En effet, on a l'identité

$$T^{\otimes(n+1)} - 1 = (T - 1) \otimes T^{\otimes n} + 1 \otimes (T - 1) \otimes T^{\otimes(n-1)} + ... + 1 \otimes 1 \otimes ... \otimes (T - 1)$$

L'opérateur d'homotopie cherché $K_{n+1}$ est alors

$$K \otimes T^{\otimes n} + 1 \otimes K \otimes T^{\otimes(n-1)} + ...\ 1 \otimes 1 \otimes ... \otimes 1 \otimes K$$

Sur cette formule, on voit que $K_{n+1}$ conserve les produits tensoriels de fonctions de Dirac et de Heaviside dans l'algèbre $\mathcal{D}^*(x_0, ..., x_n) = \mathcal{D}^*(x_0) \otimes ... \otimes \mathcal{D}^*(x_n)$ ; il induit donc un opérateur d'homotopie sur l'algèbre quotient $\mathcal{D}^*(\Delta_n)$ (cf. 2.8). Par ailleurs, les opérateurs face consistent à poser une variable $x_i$ égale à $- \infty$ et ils sont aussi compatibles avec les opérateurs $K_n$ et $K_{n+1}$. Par contre ils ne sont pas compatibles avec les opérateurs de dégénérescence.

**8.9.** L'opérateur de translation T se transporte évidemment à $\mathcal{D}^*(X)$ si X est un complexe simplicial fini comme dans le § 4. Puisque les opérateurs d'homotopie précédents $K_*$ sont compatibles avec les opérateurs face, il en résulte qu'ils définissent un opérateur d'homotopie global K sur $\mathcal{D}^*(X)$ tel que

$$dK + Kd = T - 1.$$

**8.10. THEOREME.** *Soient* X *et* Y *deux complexes* ***finis*** *et soit* z *un élément de* $\mathcal{D}^*(X) \otimes \mathcal{D}^*(Y)$. *Il existe alors un entier* m *tel que pour tout* $n \geq m$, ou $n < - m$, *on ait* $(T^n \otimes 1)(z) \in \mathcal{D}^*(X)\ \overline{\otimes}\ \mathcal{D}^*(Y)$.

*Démonstration.* On peut évidemment supposer que z s'écrit $x_1 \otimes x_2$, où

$$x_1 \in \mathcal{D}^*(X) = \mathrm{Hom}(X_\# , \mathcal{D}^*(\Delta_\#))$$

$$x_2 \in \mathcal{D}^*(Y) = \mathrm{Hom}(Y_\# , \mathcal{D}^*(\Delta_\#)$$

(Ici le foncteur Hom doit être compris comme k-module des morphismes entre ensembles **semi-simpliciaux**).

On remarque alors que l'opérateur T translate les supports singuliers des fonctions de Heaviside et de Dirac (nécessairement en nombre fini) dans les expressions figurant dans $\mathcal{D}^*(\Delta_\#)$. Si m est assez grand en valeur absolue, les supports singuliers des fonctions intervenant dans $T^m(x_1))$ et $\mathcal{D}^*(x_2)$ seront donc situés dans des intervalles disjoints, ce qui démontre le théorème.

**8.11.** Nous allons terminer ce paragraphe (et l'article) par une dernière structure algébrique que nous avons évoquée dans le § 1 qui est celle d'ADGQ **spéciale** sur $\mathcal{D}^*(X)$. Celle-ci sera définie si X est un ensemble simplicial **fini**. Dans ce cas, on sait que $\mathcal{D}^*(X)$ s'identifie à l'ensemble des applications simpliciales de $X_\natural$ dans l'ADG simpliciale $\mathcal{D}^*(\Delta_\natural)$. De même, le produit tensoriel réduit $\mathcal{D}^*(X) \mathbin{\overline{\otimes}} \mathcal{D}^*(Y)$ s'identifie à l'ensemble des applications bisimpliciales de X x Y dans $\mathcal{D}^*(\Delta) \mathbin{\overline{\otimes}} \mathcal{D}^*(\Delta)$ avec des notations évidentes. Mais $\mathcal{D}^*(\Delta_m) \mathbin{\overline{\otimes}} \mathcal{D}^*(\Delta_n)$ est engendré par le produit réduit $\mathcal{D}^*(\Delta_m) \mathbin{\overline{x}} \mathcal{D}^*(\Delta_n)$ défini en 4.5. Il faut voir maintenant que $\mathcal{D}^*(X) \mathbin{\overline{\otimes}} \mathcal{D}^*(Y)$ est engendré par le produit réduit $\mathcal{D}^*(X) \mathbin{\overline{x}} \mathcal{D}^*(Y)$, défini comme l'ensemble des **applications** bisimpliciales de X x Y dans $\mathcal{D}^*(\Delta) \mathbin{\overline{x}} \mathcal{D}^*(\Delta)$. Pour cela, nous allons raisonner par une double récurrence sur le nombre de cellules de X et Y. Pour fixer les idées, supposons que X soit obtenu à partir de X' par adjonction d'une cellule de dimension n. On a donc un diagramme cocartésien

$$\begin{array}{ccc} \partial\Delta_n & \longrightarrow & X' \\ \downarrow & & \downarrow \\ \Delta_n & \longrightarrow & X \end{array}$$

d'où on déduit un diagramme cartésien

$$\begin{array}{ccc} \mathcal{D}^*(X) \mathbin{\overline{\otimes}} \mathcal{D}^*(Y) & \longrightarrow & \mathcal{D}^*(\Delta_n) \mathbin{\overline{\otimes}} \mathcal{D}^*(Y) \\ \downarrow & & \downarrow \\ \mathcal{D}^*(X') \mathbin{\overline{\otimes}} \mathcal{D}^*(Y) & \longrightarrow & \mathcal{D}^*(\partial\Delta_n) \mathbin{\overline{\otimes}} \mathcal{D}^*(Y) \end{array}$$

Consérons maintenant un élément ω du produit tensoriel réduit $\mathcal{D}^*(X) \mathbin{\overline{\otimes}} \mathcal{D}^*(Y)$. Il s'agit de montrer qu'il est engendré par $\mathcal{D}^*(X) \mathbin{\overline{x}} \mathcal{D}^*(Y)$, qui se situe aussi dans un diagramme cartésien **d'ensembles**

$$\begin{array}{ccc}
\mathcal{D}^*(X)\ \bar{x}\ \mathcal{D}^*(Y) & \longrightarrow & \mathcal{D}^*(\Delta_n)\ \bar{x}\ \mathcal{D}^*(Y) \\
\downarrow & & \downarrow \\
\mathcal{D}^*(X')\ \bar{x}\ \mathcal{D}^*(Y) & \longrightarrow & \mathcal{D}^*(\partial\Delta_n)\ \bar{x}\ \mathcal{D}^*(Y)
\end{array}$$

D'après l'hypothèse de récurrence, chacun des trois produits tensoriels réduits n'impliquant pas X est engendré par les produits réduits correspondants. Par le lemme d'extension vu en 5.5, il suffit de considérer un élément $\sum \lambda_i f_i \otimes g_i$ du produit tensoriel réduit $\mathcal{D}^*(\Delta_n)\ \overline{\otimes}\ \mathcal{D}^*(Y)$ avec $(f_i, g_i) \in \mathcal{D}^*(\Delta_n)\ \bar{x}\ \mathcal{D}^*(Y)$ dont l'image $\sum \lambda_i \bar{f}_i \otimes g_i$ dans $\mathcal{D}^*(\partial\Delta_n)\ \overline{\otimes}\ \mathcal{D}^*(Y)$ est égale à 0. En multipliant $\bar{f}_i$ par une fonction de Heaviside adéquate (cf. 5.3), on peut même supposer que $\bar{f}_i = 0$. Le couple $(f_i, g_i)$ définit alors un élément du produit réduit $\mathcal{D}^*(X)\ \bar{x}\ \mathcal{D}^*(Y)$ de manière évidente. Ces réductions successives montrent bien que l'élément ω est engendré par $\mathcal{D}^*(X)\ \bar{x}\ \mathcal{D}^*(Y)$.